# LOCAL CENTRAL LIMIT THEOREMS, THE HIGH-ORDER CORRELATIONS OF REJECTIVE SAMPLING AND LOGISTIC LIKELIHOOD ASYMPTOTICS

By Richard Arratia, Larry Goldstein[1] and Bryan Langholz[1]

*University of Southern California*

Let $I_1, \ldots, I_n$ be independent but not necessarily identically distributed Bernoulli random variables, and let $X_n = \sum_{j=1}^n I_j$. For $\nu$ in a bounded region, a local central limit theorem expansion of $\mathbb{P}(X_n = \mathbb{E} X_n + \nu)$ is developed to any given degree. By conditioning, this expansion provides information on the high-order correlation structure of dependent, weighted sampling schemes of a population $E$ (a special case of which is simple random sampling), where a set $\mathbf{d} \subset E$ is sampled with probability proportional to $\prod_{A \in \mathbf{d}} x_A$, where $x_A$ are positive weights associated with individuals $A \in E$. These results are used to determine the asymptotic information, and demonstrate the consistency and asymptotic normality of the conditional and unconditional logistic likelihood estimator for unmatched case-control study designs in which sets of controls of the same size are sampled with equal probability.

**1. Introduction.** The unmatched case-control study is one of the most widely used designs in chronic disease epidemiologic research. Typically, a large number of individuals, the cohort or *study base*, will be observed for occurrence of a binary disease outcome. Because the number of subjects is large and only a small proportion will be *cases* that contract the disease of interest, nondiseased *controls* are sampled to serve as a comparison group. Exposure and other covariate information is then obtained for the case-control study subjects for use in statistical analyses. As an example, in a study to assess the association of a variety of hypertensive drugs and the risk of myocardio-infarction (MI), 623 MI cases who used antihypertensive drugs were identified within an HMO in Washington State. The cases were grouped by sex, 10-year age, and calendar year of MI [Psaty et al. (1995)]. For each

---

Received November 2001; revised March 2004.

[1]Supported in part by NIH Grant CA42949.

group, a number of controls from the antihypertensive drug users were sampled in a fixed proportion to the number of cases. For each case-control study member, the types of antihypertensive drugs used were ascertained through computerized records, chart review and interview. The primary method of analysis was unconditional logistic regression. It was found that risk of MI was 60% higher among calcium channel blocker users compared to either diuretics alone or compared to $\beta$-blockers, a finding that has resulted in a change in treatment strategy.

The structure of these data is *prospective* in that disease occurrence is conditional on the covariate information, and controls are randomly sampled from the pool of nondiseased. This is the structure of a nested case-control study from the study base [Mantel (1973)], which we call the *nested case-control* data model. Another way to view case-control data is *retrospectively* in which the case and control covariate values are taken to be independent realizations from their respective distributions [e.g., Breslow and Powers (1978), Prentice and Pyke (1979), Weinberg and Wacholder (1993) and Carroll, Wang and Wang (1995)]. Although the nested case-control model is used in modern texts on case-control studies in epidemiologic research [e.g., Breslow and Day (1980), Kelsey, Whittemore, Evans and Thompson (1996) and Rothman and Greenland (1998)], it has been the retrospective model that is invoked when developing estimators and analyzing their properties. However, the assumption that the case and control covariates are independent random replicates may not hold in practice. For instance, if the distribution of drug types changed during the antihypertensive drug-MI study, differences in treatment within the case and control populations would make the modeling of the covariates by a common distribution within each group untenable, so the conditions required by the retrospective model analysis would not be met. But, it seems evident that valid results can still be drawn from such a study since the assignment of drug type to subjects should not influence the association between the drug type and disease.

In this paper we develop the theory necessary to determine the asymptotic behavior of estimators of the odds ratio in the nested case-control model under general conditions on the covariates and sampling methods. We then apply this theory to the maximum conditional and unconditional logistic likelihood estimators. Although the conditional logistic likelihood gives rise to valid estimators in a wider range of case-control study settings than the unconditional (e.g., individually matched case-control designs), its asymptotic properties for "large strata" have not been studied. The path of our analysis leads us through some unexpectedly broad territory, including a high-order local central limit theorem for the Poisson–Binomial distribution and expansions for the inclusion probabilities and correlation structure of rejective sampling.

After formally introducing the problem of analysis of case-control data in Section 1.1, in Section 2 we prove Theorem 2.1, a high-order local central limit theorem for the sum $X_n$ of independent but not necessarily identically distributed Bernoulli random variables having success probability $p_j, j = 1, 2, \ldots$. This result gives an expansion to any desired order for the probability that the sum $X_n$ deviates from its mean $\mathbb{E}X_n$ by the value $\nu$, uniformly for $\nu$ in any bounded region. This result is of independent interest, as it provides a means to approximate, with rates, the Poisson–Binomial distribution, for which no simple expression exists.

In Section 3 we extend Theorem 2.1 by showing that this local central limit theorem expansion holds for the sums $X_E$ of independent Bernoulli variables with success probability

$$p_{A,\lambda} = \frac{\lambda x_A}{1 + \lambda x_A}, \qquad A \in E,$$

uniformly for all $\lambda$ in an interval bounded away from zero and infinity, under asymptotic stability conditions on the weights $x_A, A \in E$. For any $\lambda > 0$, conditioning the Bernoulli variables on the event $X_E = \eta$ gives Hájek's rejective sampling scheme $\mathbb{E}_{E,\eta}$ on $E$, where a set $\mathbf{d} \subset E$ of size $\eta$ is sampled with probability proportional to $x_\mathbf{d}$, the product of the weights $x_A$ over $A \in \mathbf{d}$. Choosing $\lambda$ so that the expected number of successes $\mathbb{E}X_E$ equals $\eta$ allows for the application of local central limit Theorem 2.1, yielding Theorem 3.1, which gives an expansion for the inclusion probabilities under the rejective sampling scheme $\mathbb{E}_{E,\eta}$. This expansion is applied in Section 4 to derive Theorem 4.1, yielding the high-order correlation structure of rejective sampling.

In Section 5 we apply the rejective sampling results to the asymptotics of estimators under the nested case-control model. Theorems 5.1 and 5.2 give the asymptotic information and demonstrate the consistency and asymptotic normality of the conditional and unconditional logistic maximum likelihood estimators, respectively.

Finally, in Section 6, we compare our approach to others, and, in particular, to the derivation of asymptotics by Prentice and Pyke (1979) under the retrospective model. Lastly, we discuss efficiency issues, extensions and directions for further research.

1.1. *The statistical model and likelihood.* The prospective logistic model for disease occurrence is as follows: with covariate vector $\mathbf{z} \in \mathbb{R}^p$, the probability of disease is

(1) $$p_\lambda(\mathbf{z}; \boldsymbol{\beta}) = \frac{\lambda x(\mathbf{z}; \boldsymbol{\beta})}{1 + \lambda x(\mathbf{z}; \boldsymbol{\beta})},$$

where $x(\mathbf{z}, \mathbf{0}) = x(\mathbf{0}, \boldsymbol{\beta}) = 1$, for all $\mathbf{z} \in \mathbb{R}^p$ and $\boldsymbol{\beta}$ in the parameter space $\mathcal{B} \subset \mathbb{R}^p$ [e.g., Breslow and Day (1980) and Cox and Snell (1989)]. The parameter $\lambda > 0$ is therefore the baseline odds and $x(\mathbf{z}, \boldsymbol{\beta})$ is the odds ratio associated with $\mathbf{z}$. The odds ratio parameter $\boldsymbol{\beta}$ is typically of primary interest.

We consider a "study base" $\mathcal{R} = \{1, \ldots, N\}$ of $N$ individuals with covariates $\mathbf{z}_j$, $j \in \mathcal{R}$, and independent failure indicators $I_j$ having marginal distribution given by (1) for some $(\lambda_0, \boldsymbol{\beta}_0)$, that is, $\mathbb{P}_{\lambda_0, \beta_0}(I_j = 1) = p_{\lambda_0}(\mathbf{z}_j; \boldsymbol{\beta}_0)$. Define $x_j(\boldsymbol{\beta}) = x(\mathbf{z}_j; \boldsymbol{\beta})$, $p_{j,\lambda}(\boldsymbol{\beta}) = 1 - q_{j,\lambda}(\boldsymbol{\beta}) = p_\lambda(\mathbf{z}_j; \boldsymbol{\beta})$; we may further suppress $\boldsymbol{\beta}_0$ and write, for example, $x_j = x_j(\boldsymbol{\beta}_0)$ and $p_{j,\lambda} = p_{j,\lambda}(\boldsymbol{\beta}_0)$. Denoting the set of indices of diseased subjects by $\mathbf{D}$, for $\mathbf{d} \subset \mathcal{R}$, the probability of observing $\mathbf{D} = \mathbf{d}$ is therefore

$$\mathbb{P}_{\lambda_0, \beta_0}(\mathbf{D} = \mathbf{d}) = \prod_{j \in \mathcal{R}} p_{j,\lambda_0}^{I_j} q_{j,\lambda_0}^{1-I_j} = \lambda_0^{|\mathbf{d}|} \, x_\mathbf{d} q_\mathcal{R},$$

where for any $F \subset \mathcal{R}$,

(2) $\quad q_F = q_F(\lambda_0, \boldsymbol{\beta}_0) = \prod_{j \in F}(1 + \lambda_0 x_j)^{-1} \quad \text{and} \quad x_\mathbf{d} = \prod_{j \in \mathbf{d}} x_j.$

When covariate values for all study base subjects are available, estimation of the unknown $\boldsymbol{\beta}_0$ (and $\lambda_0$) can be achieved by maximizing the likelihood $\mathbb{P}_{\lambda,\beta}(\mathbf{D})$. But when the study base is large or the collection of the full set of covariate values is expensive or impractical, it is natural to sample subjects to form a sampled study base $E \subset \mathcal{R}$ and use the collected covariates in the sample for the estimation of parameters. Generally, a sampling design is specified by $\pi(\mathbf{s}|\mathbf{d})$, the probability of choosing $\mathbf{s}$ as the sampled risk set $E$ when $\mathbf{d}$ is the observed set of diseased subjects.

For the calculation of a likelihood, additional information that $\mathbf{D} \in \mathcal{S}$ for some $\mathcal{S}$ may be included. Conditioning on $E$ and $\mathcal{S}$ leads to the probability

(3) $\quad \mathbb{P}_{\lambda_0, \beta_0}(\mathbf{D}|E, \mathcal{S}) = \dfrac{\lambda_0^{|\mathbf{D}|} x_\mathbf{D}(\boldsymbol{\beta}_0) \pi(E|\mathbf{D})}{\sum_{\mathbf{u} \subset \mathcal{S}} \lambda_0^{|\mathbf{u}|} x_\mathbf{u}(\boldsymbol{\beta}_0) \pi(E|\mathbf{u})}.$

A likelihood is formed by allowing the parameters in (3) to vary to obtain the likelihood function $L_{E,\mathcal{S}}(\lambda, \boldsymbol{\beta}) = \mathbb{P}_{\lambda,\beta}(\mathbf{D}|E, \mathcal{S})$.

Of particular interest for epidemiologic unmatched case-control studies is the likelihood which results from (3) when conditioning on the number of cases in the case-control set. In practice, in unmatched case-control studies one typically has information on all cases and a set of controls obtained using sampling schemes such as frequency matching, fixed size sampling, Bernoulli trials and case-base sampling [e.g., Kupper, McMichael and Spirtas (1975), Breslow and Day (1980), Wacholder, Silverman, McLaughlin and Mandel (1992) and Langholz and Goldstein (2001)]. For each of these designs, the probability $\pi(\mathbf{s}|\mathbf{d})$ is zero unless $\mathbf{s}$ contains $\mathbf{d}$, and is otherwise constant in $|\mathbf{d}|$. Then, setting $\mathcal{S} = \{\mathbf{u} \subset E : |\mathbf{u}| = \eta\}$, where $\eta = |\mathbf{D}|$, $\lambda_0$ and the sampling probabilities $\pi(\mathbf{s}|\mathbf{d})$ cancel from (3), and noting the dependence of the resulting

probability on $E$ and $\eta$ only, we define

(4) $$\mathbb{P}_{E,\eta}(\mathbf{D}) = \mathbb{P}_{\beta_0}(\mathbf{D}|E,\eta) = \frac{x_\mathbf{D}}{\sum_{\mathbf{u} \subset E\,:\,|\mathbf{u}|=\eta} x_\mathbf{u}}.$$

This is the basis for the "standard" conditional logistic likelihood $L_{E,\eta}(\boldsymbol{\beta}) = \mathbb{P}_\beta(\mathbf{D}|E,\eta)$ [e.g., Cox (1972) and Cox and Snell (1989)] for the designs mentioned above, which have log likelihood

$$\mathcal{L}_{E,\eta}(\boldsymbol{\beta}) = \sum_{A \in \mathbf{D}} \log x_A(\boldsymbol{\beta}) - \log\left\{\sum_{\mathbf{u} \subset E, |\mathbf{u}|=\eta} x_\mathbf{u}(\boldsymbol{\beta})\right\}.$$

The conditional logistic likelihood estimator $\hat{\boldsymbol{\beta}}_N$ is a value maximizing $\mathcal{L}_{E,\eta}(\boldsymbol{\beta})$.

Differentiation of an array $\mathbf{F}_\beta = \{F_{i_1,\ldots,i_a}(\boldsymbol{\beta})\} \in \mathbb{R}^{n_1 \times \cdots \times n_a}$ with respect to $\boldsymbol{\beta}$ will be denoted by "$'$", resulting in the array $\mathbf{F}'_\beta = \{F'_{j,i_1,\ldots,i_a}(\boldsymbol{\beta})\} = \{(\partial/\partial\beta_j)F_{i_1,\ldots,i_a}(\boldsymbol{\beta})\} \in \mathbb{R}^{p \times n_1 \times \cdots \times n_a}$. For $\mathbf{U} \in \mathbb{R}^{n_1 \times \cdots \times n_a}$ and $\mathbf{V} \in \mathbb{R}^{m_1 \times \cdots \times m_b}$, the tensor product $\mathbf{U} \otimes \mathbf{V} \in \mathbb{R}^{n_1 \times \cdots \times n_a \times m_1 \times \cdots \times m_b}$ has components $U_{i_1,\ldots,i_a} \cdot V_{j_1,\ldots,j_b}$, and we set $|\mathbf{U}| = \sum_{i_1,\ldots,i_a} |U_{i_1,\ldots,i_a}|$, the $L^1$ norm.

CONDITION 1.1. The real valued function $x$ is positive, three times differentiable in $\boldsymbol{\beta}$ and $0 < \inf_{|\mathbf{z}| \le c} x(\boldsymbol{\beta}_0, \mathbf{z}) \le \sup_{|\mathbf{z}| \le c} x(\boldsymbol{\beta}_0, \mathbf{z}) < \infty$ for all $c > 0$.

Under Condition 1.1, following Barlow and Prentice (1988), define the "effective covariates" $\mathbf{z}_j$ by

$$\mathbf{z}_j = \mathbf{x}'_j x_j^{-1} \in \mathbb{R}^p;$$

in the model where $x_j(\boldsymbol{\beta}) = \exp(\boldsymbol{\beta}^\mathsf{T} \mathbf{z}_j)$, we have $\mathbf{z}_j = \mathbf{z}_j$. Now for $\mathbf{u} \subset E$, define the inclusion probabilities

$$p_\mathbf{u}(\boldsymbol{\beta}) = \mathbb{P}_\beta(\mathbf{u} \subset \mathbf{D}|E,\eta) = \sum_{\mathbf{s} \supset \mathbf{u}} \mathbb{P}_\beta(\mathbf{s}|E,\eta),$$

and the inclusion probability for an individual $A$ as $p_A(\boldsymbol{\beta}) = p_{\{A\}}(\boldsymbol{\beta})$. With $I_A$ the failure indicator for $A \in E$ and suppressing the dependence of $\mathbf{Z}_A$ and $p_A$ on $\boldsymbol{\beta}$, the score $\partial \mathcal{L}_{E,\eta}(\boldsymbol{\beta})/\partial \boldsymbol{\beta}$ equals

$$\mathcal{U}_{E,\eta}(\boldsymbol{\beta}) = \sum_{A \in E} \mathbf{Z}_A(I_A - p_A) = \sum_{A \in \mathbf{D}} \mathbf{Z}_A - \mathbb{E}_\beta\left(\sum_{A \in \mathbf{D}} \mathbf{Z}_A \bigg| E, \eta\right),$$

where $\mathbb{E}_\beta$ is the expectation under $\mathbb{P}_{E,\eta}(\boldsymbol{\beta})$. Using that for a function $\mathbf{F}_\beta(\mathbf{D}) = \sum_{A \in \mathbf{D}} \mathbf{Z}_A$ we have

(5) $$\frac{\partial}{\partial \boldsymbol{\beta}} \mathbb{E}_\beta(\mathbf{F}_\beta(\mathbf{D})|E,\eta) = \mathbb{E}_\beta(\mathbf{F}_\beta(\mathbf{D})'|E,\eta) + \mathbb{E}_\beta(\mathcal{U}(\boldsymbol{\beta}) \otimes \mathbf{F}_\beta(\mathbf{D})|E,\eta),$$

with $p_A + q_A = 1$, the information $-\partial \mathcal{U}_{E,\eta}(\boldsymbol{\beta})/\partial \boldsymbol{\beta}$ is given by

$$\text{(6)} \quad \mathcal{I}_{E,\eta}(\boldsymbol{\beta}) = \sum_{A \in E} \mathbf{Z}_A^{\otimes 2} p_A q_A + \sum_{A,B \in E, A \neq B} \mathbf{Z}_A \mathbf{Z}_B^{\mathsf{T}} (p_{AB} - p_A p_B)$$

$$\text{(7)} \quad - \left( \sum_{A \in \mathbf{D}} \mathbf{Z}'_A - \sum_{A \in E} \mathbf{Z}'_A p_A \right).$$

Note that (6) contains $p_{AB} - p_A p_B$, the correlation of the joint inclusion of $A$ and $B$.

In general, we have

$$\mathbb{E}_{\beta_0} \mathcal{U}_{E,\eta}(\boldsymbol{\beta}_0) = 0,$$

since the score is the difference between a quantity and its conditional expectation. For this same reason, when taking expectation in (7), we find that

$$\mathbb{E}_{\beta_0} \mathcal{I}_{E,\eta}(\boldsymbol{\beta}_0) = \mathrm{Var}_{\beta_0} \{\mathcal{U}_{E,\eta}(\boldsymbol{\beta}_0)\}.$$

The standard likelihood argument to show the consistency of $\hat{\boldsymbol{\beta}}_N$ requires that the information $|E|^{-1} \mathcal{I}(\boldsymbol{\beta}_0)$ converge in probability. Since the information is a double sum over $A$ and $B$, the inclusion correlations $p_{AB} - p_A p_B$ need to decay at rate $|E|^{-1}$. Further, the remainder term in the Taylor expansion of the log likelihood, which is required to stay bounded in probability, contains a triple sum of terms multiplied by the third-order correlation,

$$\mathbb{E}_\beta[(I_A - p_A)(I_B - p_B)(I_C - p_C)|E, \eta];$$

hence, to satisfy the boundedness condition, such triple correlations need to decay as $|E|^{-2}$. The dependence in $\mathbb{P}_{E,\eta}$ created by having the probability of a set proportional to the product of its individual weights has been explored only under very restrictive situations [Harkness (1965) and Farewell (1979); see also Hájek (1964)]. Theorem 4.1 gives information on the rate of decay on all correlation orders, and, in particular, provides that the third-order correlation decays at the required rate. This result allows for the full treatment of the asymptotic theory for the conditional logistic maximum likelihood estimator for a large class of case-control sampling designs (Section 5).

More commonly used in practice, and making use of the same case-control subject data, is the estimator of $\boldsymbol{\beta}_0$ based on maximizing the "unconditional logistic likelihood" which, with $p_{A,\lambda}(\boldsymbol{\beta})$ as in (1) and $q_E(\lambda, \boldsymbol{\beta})$ as in (2), is given by

$$\text{(8)} \quad \tilde{L}_E(\lambda, \boldsymbol{\beta}) = \prod_{A \in E} p_{A,\lambda}(\boldsymbol{\beta})^{I_A} q_{A,\lambda}(\boldsymbol{\beta})^{1-I_A} = \lambda^{|\mathbf{D}|} x_{\mathbf{D}}(\boldsymbol{\beta}) q_E(\lambda, \boldsymbol{\beta}).$$

The unconditional logistic likelihood estimator $\tilde{\boldsymbol{\beta}}_N$ is a value maximizing $\tilde{L}_E(\boldsymbol{\beta})$. Note that, in general, $\tilde{L}_E$ is not a true likelihood when data is collected using sampling methods such as frequency matching, since the contributions from individual subjects are not independent. The asymptotic analysis of the unconditional logistic estimator is carried out in Section 5.

1.2. *The probabilistic setup.* For any set $E$ and $0 \leq \eta \leq |E|$, consider the probability measure $\mathbb{P}_{E,\eta}(\mathbf{d})$ given by (4), supported on the size $\eta$ subsets of $E$. With $I_A = \mathbf{1}(A \in \mathbf{D})$, the indicator that $A$ is included in $\mathbf{D}$, $p_A = \mathbb{E}_{E,\eta}(I_A)$, and $H \subset E$, we study high-order correlations of the form

$$(9) \qquad \mathrm{Corr}(H) = \mathbb{E}_{E,\eta}\bigg(\prod_{A \in H}(I_A - p_A)\bigg).$$

When $H = \{A, B\}$, a set of size 2, $\mathrm{Corr}(H) = p_{AB} - p_A p_B$, the covariance between the Bernoulli variables $I_A$ and $I_B$.

When $x_A = 1$ for all $A \in E$ [corresponding to $\boldsymbol{\beta}_0 = 0$ in (1)], $\mathbb{P}_{E,\eta}$ reduces to simple random sampling. In this case, when there exists $\tau \in (0, 1/2]$ such that the sampling fraction $\eta/|E| \in [\tau, 1-\tau]$, then as $|E| \to \infty$,

$$(10) \qquad \mathbb{E}_{E,\eta}[(I_A - p_A)(I_B - p_B)] = \frac{-\eta(|E|-\eta)}{|E|^2(|E|-1)} = O_\tau(|E|^{-1})$$

and

$$(11) \qquad \mathbb{E}_{E,\eta}[(I_A - p_A)(I_B - p_B)(I_C - p_C)] = O_\tau(|E|^{-2}).$$

Hence, simple random sampling has the rates needed for the stability of the information and the control on the remainder in our likelihood analysis; the exact meaning of $O_\tau$ is given in Definition 2.1.

For simple random sampling a straightforward calculation shows that

$$\mathrm{Corr}(H) = \sum_{j=0}^{|H|} \frac{\binom{|H|}{j}\binom{|E|-j}{\eta-j}}{\binom{|E|}{\eta}}\left(\frac{-\eta}{|E|}\right)^{|H|-j}.$$

Since here the weights $x_A$ are equal, we may write $\mathrm{Corr}(k)$ for the common value of $\mathrm{Corr}(H)$ for all $H$ of size $k$, and have verified for $k \leq 10$, as $|E|, \eta \to \infty$, with $\eta/|E| \to f \in (0,1)$, for $\mathcal{N}$ a standard normal variate,

$$\lim_{|E| \to \infty} |E|^{k/2}\mathrm{Corr}(k) = \mathbb{E}\mathcal{N}^k(f(f-1))^{k/2} \qquad \text{for } k \text{ even}$$

and

$$\lim_{|E| \to \infty} |E|^{(k+1)/2}\mathrm{Corr}(k)$$
$$= \tfrac{1}{3}(k-1)\mathbb{E}\mathcal{N}^{k+1}(f(f-1))^{(k-1)/2}(2f-1) \qquad \text{for } k \text{ odd}.$$

In particular, for simple random sampling we have

$$(12) \qquad \mathrm{Corr}(H) = O_{|H|,\tau}(|E|^{-(|H|+|H| \bmod 2)/2}),$$

with (10) and (11) as special cases. Theorem 4.1 shows that the orders in (12) are obtained quite generally for the weighted sampling scheme $\mathbb{P}_{E,\eta}$.

1.3. *Rejective sampling.* The scheme corresponding to the probability measure $\mathbb{P}_{E,\eta}$ is known as rejective sampling [Hájek (1964)], and as seen in Section 1.2 includes simple random sampling as a particular case. Though simple random sampling is the most ubiquitous of all statistical methods, in some cases it is not possible to take a simple random sample. For example, the inclusion of the population member $A$ might be influenced by a certain nonnegative "size" $x_A$ associated with the item $A$, where the larger the size of an item, the easier it is to locate and the higher the probability of its inclusion.

The term rejective sampling arises since $\mathbb{P}_{E,\eta}$ may be achieved by sampling $\eta$ individuals independently with replacement and rejecting those samples in which the $\eta$ individuals are not distinct. Hájek (1964) considers the inclusion probabilities, second-order correlations and asymptotic normality of sums obtained by rejective sampling.

Schemes where objects are sequentially sampled proportional to their size have been extensively studied [e.g., Rosén (1972) and Gordon (1983)]. However, rejective sampling differs from sampling sequentially proportional to size when $\eta \geq 2$, as can be seen by comparison of the general probability that a sample of size $\eta = 2$ results in the units $A$ and $B$. However, both schemes reduce to simple random sampling when the weights are constant.

**2. A high-order local central limit theorem.** The main result of this section is Theorem 2.1, a local central limit theorem expansion for the distribution of $X_n$, the sum of independent but not necessarily identically distributed indicator random variables. The first step, Lemma 2.1, is to obtain an expression for the characteristic function of the centered sum, $X_n - \mathbb{E}X_n$. In the following, we write $\Omega$ for a complex number, not necessarily the same at each occurrence, such that $|\Omega| \leq 1$.

LEMMA 2.1. *Let*

$$X_n = \sum_{j=1}^{n} I_j,$$

*where $I_j, j = 1, \ldots, n$, are independent Bernoulli variables with $\mathbb{E}I_j = p_j = 1 - q_j$; let*

(13) $$v_n^2 = \sum_{j=1}^{n} p_j q_j \quad and \quad w_n = \sum_{j=1}^{n} p_j q_j (p_j - q_j).$$

*Then, denoting the characteristic function of $X_n - \mathbb{E}X_n$ by $\phi_n(t)$, for all $n = 1, 2, \ldots$ and $|t| \leq 1$,*

(14) $$\phi_n(t) = \exp\left(-\frac{t^2 v_n^2}{2} + i\frac{t^3 w_n}{6} + \frac{t^4}{10}n\Omega\right).$$

*Furthermore, for all $t \in [-\pi, \pi]$,*

(15) $$|\phi_n(t)| \leq \exp(-t^2 v_n^2/6).$$

PROOF. The characteristic function of an indicator $I$ which has been centered by subtraction of its mean $p$ is

$$\mathbb{E}e^{it(I-p)} = e^{-itp}(q + pe^{it}) = qe^{-itp} + pe^{itq}.$$

We have for all $t$,

$$qe^{-itp} = q\left(1 - itp - \frac{t^2}{2}p^2 + i\frac{t^3}{6}p^3 + \frac{t^4}{24}p^4\Omega\right),$$

and adding the analogous expansion for $pe^{itq}$, we obtain

$$\mathbb{E}e^{it(I-p)} = 1 - \frac{t^2}{2}pq + i\frac{t^3}{6}pq(p-q) + \Omega\frac{t^4}{24}.$$

Using that $pq(p-q) \leq \sqrt{3}/18 \leq 1/9$, we have for $|t| \leq 1$,

$$\left|-\frac{t^2}{2}pq + i\frac{t^3}{6}pq(p-q) + \Omega\frac{1}{24}t^4\right| \leq \frac{t^2}{2}\frac{1}{4} + \frac{|t|^3}{6}\frac{1}{9} + \frac{t^4}{24} \leq \frac{5}{27}t^2 \leq \frac{1}{2}.$$

Applying the estimate

$$\log(1+x) = x + \Omega x^2 \qquad \forall\, |x| \leq \tfrac{1}{2},$$

we obtain that for $|t| \leq 1$,

$$\log(\mathbb{E}e^{it(I-p)}) = -\frac{t^2}{2}pq + i\frac{t^3}{6}pq(p-q) + \Omega\frac{1}{24}t^4 + \Omega\left(\frac{5}{27}t^2\right)^2$$

$$= -\frac{t^2}{2}pq + i\frac{t^3}{6}pq(p-q) + \frac{t^4}{10}\Omega,$$

and now summing,

$$\log \phi_n(t) = -\frac{t^2}{2}\sum_{j=1}^n p_j q_j + i\frac{t^3}{6}\sum_{j=1}^n p_j q_j(p_j - q_j) + n\frac{t^4}{10}\Omega.$$

Exponentiating gives (14).

To prove (15), observe that

$$|\phi_n(t)| = \prod_{j=1}^n |\mathbb{E}e^{it(I-p_j)}| = \prod_{j=1}^n (p_j^2 + q_j^2 + 2p_j q_j \cos(t))^{1/2}$$

$$= \left(\prod_{j=1}^n (1 - 2(1-\cos(t))p_j q_j)\right)^{1/2} \leq \exp\left(-(1-\cos(t))\sum_{j=1}^n p_j q_j\right),$$

and then use (13) and that $1 - \cos(t) \geq t^2/6$ for all $-\pi \leq t \leq \pi$. □

DEFINITION 2.1. For a possibly empty set of parameters $\mu$, we will write $f_n = O_\mu(g_n)$ if there exist a constant $C_\mu$ and an integer $n_\mu$, both depending only on $\mu$, such that

$$|f_n| \leq C_\mu |g_n| \quad \text{for all } n \geq n_\mu; \tag{16}$$

we write $f_n = o_\mu(g_n)$ if for every $\varepsilon > 0$, there exists $n_\mu$ such that (16) holds with $C_\mu$ replaced by $\varepsilon$. We write $f_n = \Theta_\mu(g_n)$ if $f_n = O_\mu(g_n)$ and $g_n = O_\mu(f_n)$.

In the remainder of this section, recalling $v_n^2 = \sum_{j=1}^n p_j q_j$, we will assume the following:

CONDITION 2.1. There exist $\varepsilon > 0$ and $n_\varepsilon$ such that $v_n^2 \geq \varepsilon n$ for all $n \geq n_\varepsilon$.

We will again let $\mathcal{N}$ denote a standard normal variable.

LEMMA 2.2. *Let $a_n = \sqrt{C \log n / n}$ for $C > 0$. Then under Condition 2.1,*

$$\frac{1}{2\pi} \int_{|t| \leq a_n} |t|^j \exp\left(-\frac{t^2 v_n^2}{2}\right) dt = \frac{v_n^{-(j+1)}}{\sqrt{2\pi}} (\mathbb{E}|\mathcal{N}|^j + o_{\varepsilon,j,C}(1)) \tag{17}$$

$$= \Theta_{\varepsilon,j,C}(n^{-(j+1)/2}).$$

PROOF. By the change of variable $z = v_n t$, the left-hand side of (17) becomes

$$\frac{v_n^{-(j+1)}}{\sqrt{2\pi}} \int_{|z| \leq a_n v_n} |z|^j \frac{\exp(-z^2/2)}{\sqrt{2\pi}} dz$$

$$= \frac{v_n^{-(j+1)}}{\sqrt{2\pi}} (\mathbb{E}|\mathcal{N}|^j - 2\mathbb{E}\mathcal{N}^j \mathbf{1}(\mathcal{N} > a_n v_n)),$$

but

$$\mathbb{E}\mathcal{N}^j \mathbf{1}(\mathcal{N} > a_n v_n) \leq \mathbb{E}\mathcal{N}^j \mathbf{1}(\mathcal{N} > \sqrt{C\varepsilon \log n}) = o_{\varepsilon,j,C}(1),$$

as $n \to \infty$, by the dominated convergence theorem. $\square$

For a bounded function on $[-\pi, \pi]$, define

$$\|f\|_\infty = \sup_{|t| \leq \pi} |f(t)|.$$

LEMMA 2.3. *Under Condition 2.1, for any $K > 0$ and $f(t)$ a bounded measurable function on $[-\pi, \pi]$, setting*

$$a_n = \sqrt{C \log n / n} \qquad \text{with } C \geq 6\varepsilon^{-1} K,$$

*we have*

$$\int_{a_n < |t| \leq \pi} f(t) \phi_n(t) \, dt = \|f\|_\infty O(n^{-K}).$$

PROOF. Using Lemma 2.1,

$$\left| \int_{a_n < |t| \leq \pi} f(t) \phi_n(t) \, dt \right| \leq \|f\|_\infty \int_{a_n < |t| \leq \pi} |\phi_n(t)| \, dt$$

$$\leq \|f\|_\infty \int_{a_n < |t| \leq \pi} e^{-n\varepsilon t^2/6} \, dt$$

$$\leq 2\pi \|f\|_\infty e^{-n\varepsilon a_n^2/6} \leq 2\pi \|f\|_\infty n^{-K}. \qquad \square$$

LEMMA 2.4. *Let $\phi_n(t)$ be the characteristic function of the sum of $n$ independent centered Bernoulli variables, and suppose that Condition 2.1 holds. Define for $j \geq 0$,*

$$(18) \qquad \mathcal{I}_{n,j} = \frac{1}{2\pi} \int_{|t| \leq \pi} t^j \phi_n(t) \, dt \quad and \quad \mathcal{I}_{n,j}^0 = \frac{\mathcal{I}_{n,j}}{\mathcal{I}_{n,0}}.$$

*Then for $j$ even,*

$$(19) \quad \mathcal{I}_{n,j} = \Theta_{\varepsilon,j}(n^{-(j+1)/2}) \quad and \quad \mathcal{I}_{n,j}^0 = v_n^{-j} \mathbb{E}\mathcal{N}^j + o_{\varepsilon,j}(n^{-j/2}),$$

*and for $j$ odd,*

$$\mathcal{I}_{n,j} = O_{\varepsilon,j}(n^{-(j+2)/2});$$

*in particular,*

$$\mathcal{I}_{n,j}^0 = O_{\varepsilon,j}(n^{-(j+j \bmod 2)/2}) \qquad \text{for all } j \geq 0.$$

PROOF. Let $a_n = \sqrt{C \log n / n}$ with $C = 6\varepsilon^{-1}((j+3)/2)$. Lemma 2.3 yields

$$\int_{a_n < |t| \leq \pi} t^j \phi_n(t) \, dt = \pi^j O(n^{-(j+3)/2}),$$

so it suffices to consider the region $|t| \leq a_n$. Take $n_{\varepsilon,j}$ so that for $n \geq n_{\varepsilon,j}$, $a_n \leq 1$ and $na_n^3 \leq 3$. Since $|t| \leq a_n \leq 1$, (14) of Lemma 2.1 gives

$$\phi_n(t) = \exp\left(-\frac{t^2 v_n^2}{2} + i\frac{t^3 w_n}{6} + nt^4 \frac{1}{10}\Omega\right)$$

$$= \exp\left(-\frac{t^2 v_n^2}{2}\right) \exp\left(i\frac{t^3 w_n}{6} + nt^4 \frac{1}{10}\Omega\right).$$

For $n \geq n_{\varepsilon,j}$, we see that $na_n^4/10 \leq na_n^3/6 \leq 1/2$. Therefore, for $|t| \leq a_n$, $|x| \leq 1$, where $x = it^3 w_n/6 + nt^4 \Omega/10$. Now using the fact that for $|x| \leq 1$,
$$e^x = 1 + x + O(x^2),$$
we have
$$\phi_n(t) = \exp\left(-\frac{t^2 v_n^2}{2}\right)\left(1 + i\frac{t^3 w_n}{6}\right) + \exp\left(-\frac{t^2 v_n^2}{2}\right) O(nt^4 + t^6 w_n^2).$$

Lemma 2.2 shows that the second term contributes $O_{\varepsilon,j}(n^{-(j+3)/2})$ to $\mathcal{I}_{n,j}$, since
$$n \int_{|t| \leq a_n} |t|^{j+4} \exp\left(-\frac{t^2 v_n^2}{2}\right) = n\Theta_{\varepsilon,j}(n^{-(j+5)/2}) = \Theta_{\varepsilon,j}(n^{-(j+3)/2})$$
and
$$w_n^2 \int_{|t| \leq a_n} |t|^{j+6} \exp\left(-\frac{t^2 v_n^2}{2}\right) = O_{\varepsilon,j}(n^2 n^{-(j+7)/2}) = O_{\varepsilon,j}(n^{-(j+3)/2}).$$

Now focusing on the contribution from the first term, using $v_n^2 = \Theta(n)$ by Condition 2.1, symmetry and Lemma 2.2 for $j$ even, we have
$$v_n^{j+1} \mathcal{I}_{n,j} = \frac{v_n^{j+1}}{2\pi} \int_{|t| \leq a_n} t^j \exp\left(-\frac{t^2 v_n^2}{2}\right)\left(1 + i\frac{t^3}{6} w_n\right) dt + O_{\varepsilon,j}(n^{-1})$$
$$= \frac{v_n^{j+1}}{2\pi} \int_{|t| \leq a_n} t^j \exp\left(-\frac{t^2 v_n^2}{2}\right) dt + O_{\varepsilon,j}(n^{-1})$$
$$= \frac{1}{\sqrt{2\pi}}(\mathbb{E}\mathcal{N}^j + o_{\varepsilon,j}(1)),$$
yielding (19).

For $j$ odd, again using symmetry,
$$n^{(j+2)/2} \mathcal{I}_{n,j} = n^{(j+2)/2} \frac{1}{2\pi} \int_{|t| \leq a_n} t^j \exp\left(-\frac{t^2 v_n^2}{2}\right)\left(1 + i\frac{t^3 w_n}{6}\right) dt + O_{\varepsilon,j}(n^{-1/2})$$
$$= n^{(j+2)/2} \frac{iw_n}{12\pi} \int_{|t| \leq a_n} t^{j+3} \exp\left(-\frac{t^2 v_n^2}{2}\right) dt + O_{\varepsilon,j}(n^{-1/2})$$
$$= \frac{i(w_n/n)}{6\sqrt{2\pi}} \left(\frac{n}{v_n^2}\right)^{(j+4)/2} (\mathbb{E}\mathcal{N}^{j+3} + o_{\varepsilon,j}(1)) + O_{\varepsilon,j}(n^{-1/2});$$
the right-hand side is now seen to be $O_{\varepsilon,j}(1)$. □

For $\mathbb{E}X_n + \nu$ an integer, define
$$(20) \qquad f_{n,\nu} = \mathbb{P}(X_n = \mathbb{E}X_n + \nu).$$
The following theorem gives a high-order local central limit for the probabilities of such deviations from the mean $\mathbb{E}X_n$.

THEOREM 2.1. *Let $I_1, I_2, \ldots$ be independent Bernoulli variables with $p_j = \mathbb{E} I_j$ satisfying Condition 2.1. For any nonnegative integer $s$, define*

$$\tag{21} m_\nu(s) = \sum_{j=0}^{s} \frac{(-i\nu)^j}{j!} \mathcal{I}_{n,j}.$$

*Then for given $\kappa$ and even $s$,*

$$f_{n,\nu} = m_\nu(s) + \Theta_{\varepsilon, \kappa, s}(n^{-(s+3)/2}) \qquad \text{for all } |\nu| \leq \kappa \text{ with } \mathbb{E} X_n + \nu \in \mathbf{N}.$$

PROOF. Let

$$R_{n,\nu} = f_{n,\nu} - \sum_{j=0}^{s} \frac{(-i\nu)^j}{j!} \mathcal{I}_{n,j}, \qquad g(x) = e^x - \sum_{j=0}^{s} \frac{x^j}{j!}$$

and $a_n = \sqrt{C \log n / n}$ with $C = 6\varepsilon^{-1}(s+2)/2$. By the inversion formula,

$$f_{n,\nu} = \frac{1}{2\pi} \int_{|t| \leq \pi} e^{-it\nu} \phi_n(t) \, dt,$$

so

$$2\pi R_{n,\nu} = \int_{|t| \leq \pi} g(-it\nu) \phi_n(t) \, dt$$

$$= \int_{|t| \leq a_n} g(-it\nu) \phi_n(t) \, dt + \int_{a_n < |t| \leq \pi} g(-it\nu) \phi_n(t) \, dt.$$

Since $|\phi_n(t)| \leq 1$, $|\nu| \leq \kappa$ and $|g(x)| \leq C_s |x|^{s+1}$ for $|x| \leq \pi$, the first integral is bounded by

$$\int_{|t| \leq a_n} |g(-it\nu)| \, dt \leq 2a_n \sup_{|t| \leq a_n} |g(-it\nu)|$$

$$\leq 2 C_s a_n^{s+2} \kappa^{s+1} = O_{\varepsilon, \kappa, s}\left(\left(\frac{\log(n)}{n}\right)^{(s+2)/2}\right).$$

Since $\sup_{|t| \leq \pi} |g(-i\nu t)| \leq C_s (\pi \kappa)^{s+1}$, Lemma 2.3 shows that the second integral is $O_{\kappa,s}(n^{-(s+2)/2}) \subset O_{\varepsilon, \kappa, s}((\log n / n)^{(s+2)/2})$. Consequently, for all $s$, we obtain

$$f_{n,\nu} = \sum_{j=0}^{s+2} \frac{(-i\nu)^j}{j!} \mathcal{I}_{n,j} + O_{\varepsilon, \kappa, s}\left(\left(\frac{\log n}{n}\right)^{(s+4)/2}\right).$$

When $s$ is even, we have by Lemma 2.4,

$$\mathcal{I}_{n,s+1} = O_{\varepsilon,s}(n^{-(s+3)/2}) \quad \text{and} \quad \mathcal{I}_{n,s+2} = \Theta_{\varepsilon,s}(n^{-(s+3)/2});$$

we now obtain the result by observing

$$O_{\varepsilon,\kappa,s}\left(\left(\frac{\log n}{n}\right)^{(s+4)/2}\right) \subset O_{\varepsilon,\kappa,s}(n^{-(s+3)/2}). \qquad \square$$

**3. Finite population sampling and inclusion probabilities.** To extend the results of Section 2, let there be given for all $A \in \mathbf{N} = \{0, 1, \dots\}$ a "weight" $x_A \geq 0$, and for $\lambda > 0$, let $T_\lambda$ be the measure under which $I_A$ for $A \in \mathbf{N}$ are independent indicator variables with success probability

$$p_{A,\lambda} = \frac{\lambda x_A}{1 + \lambda x_A}. \tag{22}$$

The case considered in Section 2 corresponds to $\lambda = 1$ and $x_j = p_j/q_j$.

We will assume that the $x_A$ weights are "asymptotically stable" in the following sense.

CONDITION 3.1. For all $\delta \in (0,1)$, there exist $\varepsilon \in (0,1)$ and $n \geq 1$ such that for any finite $E \subset \mathbf{N}$ with $|E| \geq n$,

$$\frac{1}{|E|} \sum_{A \in E} \mathbf{1}(x_A \in [\varepsilon, \varepsilon^{-1}]) \geq 1 - \delta. \tag{23}$$

Now let $p_{A,\lambda} + q_{A,\lambda} = 1$,

$$X_E = \sum_{A \in E} I_A, \qquad v_{E,\lambda}^2 = \sum_{A \in E} p_{A,\lambda} q_{A,\lambda}, \tag{24}$$

and with $T_\lambda(X_E)$ denoting the expectation of $X_E$ with respect to $T_\lambda$ and $T_\lambda(X_E) + \nu$ an integer, set

$$f_{E,\lambda,\nu} = T_\lambda(X_E = T_\lambda(X_E) + \nu). \tag{25}$$

In this section we will provide a local central limit theorem expansion for the probabilities in (25) which holds uniformly for $\lambda$ in an interval bounded away from zero and infinity. Conditioning $T_\lambda$ to have exactly $\eta$ successes over $E$ yields $\mathbb{E}_{E,\eta}$ (Lemma 3.5), and by selecting the $\lambda$ which yields $T_\lambda(X_E) = \eta$, we obtain a high-order expansion for the probability that $A$ is included in a sample with distribution $\mathbb{P}_{E,\eta}$.

With $f_{E,\nu}$ a real valued function defined on finite subsets $E \subset \mathbf{N}$ and $\nu \in \mathbb{R}$, for a possibly empty collection of parameters $\mu$, we say

$$f_{E,\nu} = O_\mu(g_E)$$

if there exist $C_\mu$ and $n_\mu$ such that

$$|f_{E,\nu}| \leq C_\mu |g_E| \qquad \text{for all } |E| \geq n_\mu.$$

We say $f_{E,\nu} = \Theta_\mu(g_E)$ when $f_{E,\nu} = O_\mu(g_E)$ and $g_E = O_\mu(f_{E,\nu})$. Note that if $H$ and $G$ are any fixed finite subsets of $\mathbf{N}$, then $f_{E,\nu} = O_\mu(|E|^{-a})$ if and only if $f_{E,\nu} = O_\mu(|(E \setminus H) \cup G|^{-a})$.

To see that Condition 3.1 implies Condition 2.1 in Section 2, uniformly for $\lambda$ in an interval bounded away from zero and infinity, we have:

LEMMA 3.1. *Let Condition* 3.1 *hold and* $\gamma \in (0,1]$. *Then there exist* $\varepsilon_\gamma > 0$ *and* $n_\gamma$ *such that*

$$v_{E,\lambda}^2 \geq \varepsilon_\gamma |E| \qquad \text{for all } \lambda \in [\gamma, 1/\gamma] \text{ and } |E| \geq n_\gamma.$$

PROOF. Letting $\delta \in (0,1), \varepsilon \in (0,1]$ and $n$ be any values satisfying (23), and set

$$\varepsilon_\gamma = \frac{(1-\delta)\gamma\varepsilon}{(1+\gamma\varepsilon)(1+\gamma^{-1}\varepsilon^{-1})} \quad \text{and} \quad n_\gamma = n.$$

Then for any $|E| \geq n_\gamma$ and $\lambda \in [\gamma, 1/\gamma]$,

$$v_{E,\lambda}^2 = \sum_{A \in E} \frac{\lambda x_A}{1+\lambda x_A} \frac{1}{1+\lambda x_A} \geq \sum_{A \in E} \frac{\gamma x_A}{1+\gamma x_A} \frac{1}{1+\gamma^{-1} x_A} \geq \varepsilon_\gamma |E|. \qquad \square$$

Now let $\phi_{E,\lambda}$ be the characteristic function of $X_E - T_\lambda(X_E)$ under the measure $T_\lambda$, and in parallel to (18) and (21), write

$$\mathcal{I}_{E,\lambda,j} = \frac{1}{2\pi} \int_{|t| \leq \pi} t^j \phi_{E,\lambda}(t)\, dt \quad \text{and} \quad m_{E,\lambda,\nu}(s) = \sum_{j=0}^{s} \frac{(-i\nu)^j}{j!} \mathcal{I}_{E,\lambda,j}.$$

LEMMA 3.2. *Let Condition* 3.1 *be satisfied and* $\gamma \in (0,1]$. *Then for all* $\lambda \in [\gamma, 1/\gamma]$, *for* $j$ *even*,

(26) $\quad \mathcal{I}_{E,\lambda,j} = \Theta_{\gamma,j}(|E|^{-(j+1)/2}) \quad$ *and* $\quad \mathcal{I}_{E,\lambda,j}^0 = v_{E,\lambda}^{-j} \mathbb{E} \mathcal{N}^j + o_{\gamma,j}(|E|^{-j/2}),$

*and for* $j$ *odd*,

(27) $$\mathcal{I}_{E,\lambda,j} = \Theta_{\gamma,j}(|E|^{-(j+2)/2});$$

*in particular, for all* $j$,

(28) $\quad \mathcal{I}_{E,\lambda,j}^0 \equiv \dfrac{\mathcal{I}_{E,\lambda,j}}{\mathcal{I}_{E,\lambda,0}} = O_{\gamma,j}(|E|^{-(j+j \bmod 2)/2}) \qquad$ *for* $j \geq 0$.

*Further, for given* $\kappa$ *and even* $s$, *for* $T_\lambda(X_E) + \nu \in \mathbf{N}$,

(29) $\quad f_{E,\lambda,\nu} = m_{E,\lambda,\nu}(s) + \Theta_{\gamma,\kappa,s}(|E|^{-(s+3)/2}) \qquad$ *for all* $|\nu| \leq \kappa$.

PROOF. Lemma 3.1 in conjunction with Lemma 2.4 gives (26)–(28), and in conjunction with Theorem 2.1 gives (29). $\square$

Now for $t = \{0, 1, \dots\}$ let

(30)
$$\Psi^t f_{E,\nu} = f_{E,\nu-t},$$
$$\Delta^0 f_{E,\nu} = f_{E,\nu}, \qquad \Delta f_{E,\nu} = f_{E,\nu} - f_{E,\nu-1} \quad \text{and} \quad \Delta^{t+1} = \Delta \Delta^t.$$

For $q$ a nonnegative integer, the following classes $\mathcal{G}_\mu^q$ of functions $f_{E,\nu}$ play a crucial role:

(31) $\quad \mathcal{G}_\mu^q = \{f_{E,\nu} : \forall t \geq 0, \Delta^t f_{E,\nu} = O_{\mu,t,\nu}(|E|^{-(t+q+(t+q) \bmod 2)/2})\}.$

LEMMA 3.3. *Let $p \leq q$ be nonnegative integers, and suppose that $f_{E,\nu} \in \mathcal{G}_\mu^p$ and $g_{E,\nu} \in \mathcal{G}_\mu^q$. Then*

(32) $$\mathcal{G}_\mu^p \supset \mathcal{G}_\mu^q,$$

(33) $$af_{E,\nu},\ f_{E,\nu} + g_{E,\nu} \in \mathcal{G}_\mu^p,$$

(34) $$\forall t \geq 0 \quad \Delta^t f_{E,\nu} \in \mathcal{G}_\mu^{t+p},$$

(35) $$\forall t \geq 0 \quad f_{E,\nu} - f_{E,\nu-t} \in \mathcal{G}_\mu^{p+1},$$

(36) $$\forall t \geq 0 \quad \Psi^t f_{E,\nu} \in \mathcal{G}_\mu^p,$$

(37) $$f_{E,\nu} g_{E,\nu} \in \mathcal{G}_\mu^{p+q}.$$

PROOF. Without loss of generality take $\mu = \varnothing$. Equation (32) follows since $p + j + (p+j)$ mod 2 is increasing in $p$. Equation (33) follows by (32) and the linearity of $\Delta$. Equation (34) follows from the definition of $\mathcal{G}^q$.

For (35), write
$$f_{E,\nu} - f_{E,\nu-t} = \sum_{j=0}^{t-1} f_{E,\nu-j} - f_{E,\nu-j-1} = \sum_{j=0}^{t-1} \Delta f_{E,\nu-j};$$
by (34), the summands are in $\mathcal{G}^{p+1}$, and hence by (33), so is the sum itself, proving (35). Now $\Psi f_{E,\nu} = f_{E,\nu-1} = f_{E,\nu} - \Delta f_{E,\nu} \in \mathcal{G}^p$ by (34) and (33); the case for general $t$ in (36) follows by induction.

The verification of equation (37) can be accomplished using the fact that $\Delta^t \Psi^j = \Psi^j \Delta^t$ for all nonnegative $j, t$ and the following product rule which can be easily proved by induction:
$$\Delta^t(f_{E,\nu} g_{E,\nu}) = \sum_{0 \leq j \leq t} \binom{t}{j} (\Psi^{t-j} \Delta^j f_{E,\nu})(\Delta^{t-j} g_{E,\nu}). \qquad \square$$

For notational ease, we suppress the variable $s$ in the quantity $m_{E,\lambda,\nu}^0$ defined below. Lemmas 3.3 and 3.2 have the following consequence.

LEMMA 3.4. *Let Condition 3.1 hold and $\gamma \in (0,1]$. Then for all $\lambda \in [\gamma, 1/\gamma]$,*

(38) $$m_{E,\lambda,\nu}^0 \equiv \sum_{j=0}^{s} \frac{(-i\nu)^j}{j!} \mathcal{I}_{E,\lambda,j}^0 \in \mathcal{G}_{\gamma,s}^0.$$

*Further, defining*

(39) $$n_{E,\lambda,\nu}^0 = m_{E,\lambda,\nu}^0 - p_{A,\lambda} \Delta m_{E,\lambda,\nu}^0,$$

*we have*

(40) $$n_{E,\lambda,\nu}^0 - 1 = O_{\gamma,s,\nu}(|E|^{-1})$$

*and*

(41) $$n^0_{E,\lambda,\nu} - 1 \in \mathcal{G}^0_{\gamma,s}.$$

PROOF. Note that $\Delta^t \nu^j = 0$ for $j < t$, and hence for $0 \le t \le s$,

$$\begin{aligned}
\Delta^t m^0_{E,\lambda,\nu} &= \sum_{j=t}^{s} \frac{(-i)^j \Delta^t \nu^j}{j!} \mathcal{I}^0_{E,\lambda,j} \\
&= \sum_{j=t}^{s} \frac{(-i)^j \Delta^t \nu^j}{j!} O_{\gamma,j}(|E|^{-(j+j \bmod 2)/2}) \\
&= O_{\gamma,s,t,\nu}(|E|^{-(t+t \bmod 2)/2}).
\end{aligned}$$

For $t > s$, $\Delta^t m^0_{E,\lambda,\nu} = 0$. This proves (38).

Now (40) follows from $\mathcal{I}^0_{E,\lambda,0} = 1$,

$$m^0_{E,\lambda,\nu} = 1 + O_{\gamma,s,\nu}(|E|^{-1}) \quad \text{and} \quad \Delta m^0_{E,\lambda,\nu} = O_{\gamma,s,\nu}(|E|^{-1}).$$

By (38), we have $m_{E,\lambda,\nu} \in \mathcal{G}^0_{\gamma,s}$ and $\Delta m_{E,\lambda,\nu} \in \mathcal{G}^1_{\gamma,s}$, and (32) and (33) of Lemma 3.3 give $n_{E,\lambda,\nu} \in \mathcal{G}^0_{\gamma,s}$ upon noting that $p_{A,\lambda}$ is constant in $\nu$. Since $1 \in \mathcal{G}^0_{\gamma,s}$, applying (33) again gives (41). □

Let $E$ be a finite subset of $\mathbf{N}$, and recall $x_\mathbf{d} = \prod_{A \in \mathbf{d}} x_A$ and the probability distribution $\mathbb{E}_{E,\eta}$ given in (4). For convenience, we will write, for instance, $\mathbb{P}_{E,\eta}(A)$ in place of $\mathbb{P}_{E,\eta}(A \in \mathbf{D})$, or $\mathbb{P}_{E,\eta}(\mathbf{s})$ for $\mathbb{P}_{E,\eta}(\mathbf{s} \subset \mathbf{D})$. Also recall the product measure $T_\lambda$ with marginals given by (22), such that for all $\mathbf{d} \subset E$,

(42) $$T_\lambda(\{A \in E : I_A = 1\} = \mathbf{d}) = \lambda^{|\mathbf{d}|}\left(\prod_{A \in E} \frac{1}{1 + \lambda x_A}\right) x_\mathbf{d}.$$

The following lemma provides a key relation between $T_\lambda$ and $\mathbb{P}_{E,\eta}$; the quantity $X_E$ is as in (24).

LEMMA 3.5. *For any $(E, \eta)$ with $0 \le \eta \le |E|$, $\mathbf{d} \subset E$ with $|\mathbf{d}| = \eta$ and $\lambda > 0$,*

(43) $$\mathbb{P}_{E,\eta}(\mathbf{d}) = T_\lambda(\{A \in E : I_A = 1\} = \mathbf{d} \mid X_E = \eta),$$

*and for $A \in E$ and $F = E \setminus A$,*

(44) $$\mathbb{P}_{E,\eta}(A) = \frac{p_{A,\lambda} T_\lambda(X_F = \eta - 1)}{p_{A,\lambda} T_\lambda(X_F = \eta - 1) + q_{A,\lambda} T_\lambda(X_F = \eta)}.$$

PROOF. Summing (42) over subsets of $E$ of size $\eta$ gives

$$(45) \qquad T_\lambda(X_E = \eta) = \lambda^\eta \left( \prod_{A \in E} \frac{1}{1 + \lambda x_A} \right) \sum_{\mathbf{u} \subset E, |\mathbf{u}| = \eta} x_{\mathbf{u}},$$

and now, since $|\mathbf{d}| = \eta$, division of (42) by (45) yields (43). Next,

$$\mathbb{P}_{E,\eta}(A) = T_\lambda(I_A = 1 | X_E = \eta) = \frac{T_\lambda(I_A = 1, X_E = \eta)}{T_\lambda(X_E = \eta)}$$

$$= \frac{T_\lambda(I_A = 1, X_F = \eta - 1)}{T_\lambda(I_A = 1, X_F = \eta - 1) + T_\lambda(I_A = 0, X_F = \eta)},$$

and (44) now follows using the independence of the variables $I_A$ and $X_F$ under $T_\lambda$. □

In the following, for $\tau \in (0, 1/2]$, let

$$\mathcal{E}_\tau = \{(E, \eta) : \tau \le \eta/|E| \le 1 - \tau\}.$$

LEMMA 3.6. *Suppose that Condition 3.1 is satisfied. Then for all $\tau \in (0, 1/2]$, there exist $\gamma_\tau \in (0, 1]$ and $n_\tau$ depending only on $\tau$ such that for all $(E, \eta) \in \mathcal{E}_\tau$ with $|E| \ge n_\tau$, there exists a unique solution $\boldsymbol{\lambda} = \lambda(E, \eta)$ to the equation*

$$(46) \qquad h_E(\lambda) = \frac{\eta}{|E|}, \qquad \text{where } h_E(\lambda) = \frac{1}{|E|} \sum_{A \in E} p_{A,\lambda}$$

*and*

$$(47) \qquad \lambda(E, \eta) \in [\gamma_\tau, 1/\gamma_\tau].$$

PROOF. Let $\delta = (1/2) \min\{\tau, 1 - 2\tau\}$ and take $\varepsilon$ and $n_\tau = n$ satisfying (23) for this $\delta$. Then for all $|E| \ge n_\tau$ and $\lambda > 0$,

$$(48) \qquad (1 - \delta) \frac{\lambda \varepsilon}{1 + \lambda \varepsilon} \le \frac{1}{|E|} \sum_{A \in E} \frac{\lambda x_A}{1 + \lambda x_A} \le (1 - \delta) \frac{\lambda \varepsilon^{-1}}{1 + \lambda \varepsilon^{-1}} + \delta.$$

Hence, $h_E(\lambda)$, continuous and strictly increasing on $[0, \infty)$ as a function of $\lambda$, satisfies

$$\lim_{\lambda \to 0} h_E(\lambda) \le \delta \quad \text{and} \quad \lim_{\lambda \to \infty} h_E(\lambda) \ge 1 - \delta.$$

Since $\delta < \tau \le \eta/|E| \le 1 - \tau < 1 - \delta$, there exists a unique value $\lambda(E, \eta)$ in $(0, \infty)$ for which $h_E(\lambda)$ takes on the value $\eta/|E|$.

Since $\eta/|E| \in [\tau, 1-\tau]$ and $\lambda(E,\eta)$ solves (46), by (48)

$$(1-\delta)\frac{\lambda(E,\eta)\varepsilon}{1+\lambda(E,\eta)\varepsilon} \leq 1-\tau \quad \text{and} \quad \tau \leq (1-\delta)\frac{\lambda(E,\eta)\varepsilon^{-1}}{1+\lambda(E,\eta)\varepsilon^{-1}} + \delta,$$

yielding, respectively,

$$\lambda(E,\eta) \leq \frac{1-\tau}{\varepsilon(\tau-\delta)} \quad \text{and} \quad \lambda(E,\eta) \geq \frac{\varepsilon(\tau-\delta)}{1-\tau}.$$

Verifying that $0 < (\tau-\delta)/(1-\tau) \leq 1$ completes the proof of claim (47). $\square$

THEOREM 3.1. *Suppose that Condition* 3.1 *is satisfied and let be given* $\tau \in (0, 1/2]$, $\kappa \in \mathbf{N}$ *and $s$ even. Then there exists $n_\tau$ such that for all $(E,\eta) \in \mathcal{E}_\tau$ with $|E| \geq n_\tau$, $\boldsymbol{\lambda} = \lambda(E,\eta)$ exists, and for all $|k| \leq \kappa$,*

$$(49) \quad \mathbb{P}_{E,\eta+k}(A) = p_{A,\boldsymbol{\lambda}} m^0_{F,\boldsymbol{\lambda},\nu-1} \sum_{l=0}^{s/2} (-1)^l (n^0_{F,\boldsymbol{\lambda},\nu} - 1)^l + \Theta_{\tau,\kappa,s}(|E|^{-(s+2)/2}),$$

*for all $A \in E$, where $F = E \setminus A$, $\nu = k + p_{A,\boldsymbol{\lambda}}$, and $m^0_{F,\lambda,\nu}$ and $n^0_{F,\lambda,\nu}$ are defined in* (38) *and* (39). *In particular, for all $|\nu| \leq \kappa$,*

$$(50) \quad \mathbb{P}_{E,\eta+k}(A) = p_{A,\boldsymbol{\lambda}} + O_{\tau,\kappa}(|E|^{-1}) \quad \text{and}$$

$$\Delta \mathbb{P}_{E,\eta+k}(A) = p_{A,\boldsymbol{\lambda}} q_{A,\boldsymbol{\lambda}} \mathcal{I}^0_{F,\boldsymbol{\lambda},2} + O_{\tau,\kappa}(|E|^{-2}).$$

PROOF. By Lemma 3.6, the solutions $\boldsymbol{\lambda} = \lambda(E,\eta)$ exist for all $(E,\eta) \in \mathcal{E}_\tau$ with $|E| \geq n_\tau$ and lie in an interval $[\gamma_\tau, 1/\gamma_\tau]$ for some $\gamma_\tau \in (0,1]$ depending only on $\tau$.

Hence, first applying Lemma 3.5,

$$\mathbb{P}_{E,\eta+k}(A) = \frac{p_{A,\lambda} T_\lambda(X_F = \eta+k-1)}{p_{A,\lambda} T_\lambda(X_F = \eta+k-1) + q_{A,\lambda} T_\lambda(X_F = \eta+k)}$$

(for all $\lambda > 0$)

$$= \frac{p_{A,\boldsymbol{\lambda}} T_{\boldsymbol{\lambda}}(X_F = \eta+k-1)}{p_{A,\boldsymbol{\lambda}} T_{\boldsymbol{\lambda}}(X_F = \eta+k-1) + q_{A,\boldsymbol{\lambda}} T_{\boldsymbol{\lambda}}(X_F = \eta+k)}$$

(upon setting $\lambda = \boldsymbol{\lambda}$).

Since $\eta + k = T_{\boldsymbol{\lambda}}(X_E) + k = T_{\boldsymbol{\lambda}}(T_F) + p_{A,\boldsymbol{\lambda}} + k = T_{\boldsymbol{\lambda}}(T_F) + \nu$,

$$\mathbb{P}_{E,\eta+k}(A) = \frac{p_{A,\boldsymbol{\lambda}} f_{F,\boldsymbol{\lambda},\nu-1}}{p_{A,\boldsymbol{\lambda}} f_{F,\boldsymbol{\lambda},\nu-1} + q_{A,\boldsymbol{\lambda}} f_{F,\boldsymbol{\lambda},\nu}}.$$

Letting $\Theta(s/2) = \Theta_{\tau,\kappa,s}(|E|^{-s/2})$ for short and applying Lemma 3.2, this probability equals

$$\frac{p_{A,\boldsymbol{\lambda}} m_{F,\boldsymbol{\lambda},\nu-1} + \Theta((s+3)/2)}{p_{A,\boldsymbol{\lambda}} m_{F,\boldsymbol{\lambda},\nu-1} + q_{A,\boldsymbol{\lambda}} m_{F,\boldsymbol{\lambda},\nu} + \Theta((s+3)/2)}$$

$$= \frac{p_{A,\boldsymbol{\lambda}} m^0_{F,\boldsymbol{\lambda},\nu-1} + \Theta((s+2)/2)}{p_{A,\boldsymbol{\lambda}} m^0_{F,\boldsymbol{\lambda},\nu-1} + q_{A,\boldsymbol{\lambda}} m^0_{F,\boldsymbol{\lambda},\nu} + \Theta((s+2)/2)}$$

$$[\text{since } \mathcal{I}_{F,\boldsymbol{\lambda},0} = \Theta_\tau(|E|^{-1/2})]$$

$$= \frac{p_{A,\boldsymbol{\lambda}} m^0_{F,\boldsymbol{\lambda},\nu-1} + \Theta((s+2)/2)}{n^0_{F,\boldsymbol{\lambda},\nu} + \Theta((s+2)/2)}$$

$$= \frac{p_{A,\boldsymbol{\lambda}} m^0_{F,\boldsymbol{\lambda},\nu-1}}{n^0_{F,\boldsymbol{\lambda},\nu}} + \Theta((s+2)/2)$$

$$= \frac{p_{A,\boldsymbol{\lambda}} m^0_{F,\boldsymbol{\lambda},\nu-1}}{1 + (n^0_{F,\boldsymbol{\lambda},\nu} - 1)} + \Theta((s+2)/2).$$

Equation (40) of Lemma 3.4 gives $n^0_{F,\boldsymbol{\lambda},\nu} - 1 = O_{\tau,\nu,s}(|E|^{-1})$, hence a Taylor expansion in $x$ of the quotient $1/(1+x)$ to order $s/2$ yields an error term of order $(n^0_{F,\boldsymbol{\lambda},\nu} - 1)^{s/2+1} = O_{\tau,\nu,s}(|E|^{-(s+2)/2})$, and therefore (49).

Using $s = 2$ in (49) and collecting terms of order $O_{\tau,\nu}(|E|^{-2})$, we obtain

$$\mathbb{P}_{E,\eta+k}(A) = p_{A,\boldsymbol{\lambda}}(1 + q_{A,\boldsymbol{\lambda}}(i\mathcal{I}^0_{F,\boldsymbol{\lambda},1} + (\nu - 1/2)\mathcal{I}^0_{F,\boldsymbol{\lambda},2}) + O_{\tau,\nu}(|E|^{-2})),$$

proving (50). □

Under the hypotheses of Theorem 3.1 we have the following:

COROLLARY 3.1.

(51) $$\mathbb{P}_{E,\eta+k}(A) \in \mathcal{G}^0_\tau.$$

PROOF. For $t = 0$, $\Delta^0 \mathbb{P}_{E,\eta+k}(A) = \mathbb{P}_{E,\eta+k}(A) = O_{\tau,k}(1)$. Given arbitrary $t \geq 1$, take

$$s = t + t \bmod 2 - 2 \quad \text{and} \quad \kappa = t.$$

Since $m^0_{F,\boldsymbol{\lambda},\nu} \in \mathcal{G}^0_\tau$ and (41) of Lemma 3.4 gives $n^0_{F,\boldsymbol{\lambda},\nu} - 1 \in \mathcal{G}^0_\tau$, repeated application of Lemma 3.3 shows

$$p_{A,\boldsymbol{\lambda}} m_{F,\boldsymbol{\lambda},\nu} \sum_{l=1}^{s/2} (n^0_{F,\boldsymbol{\lambda},\nu} - 1)^l \in \mathcal{G}^0_\tau.$$

For the error term, by the choice $s + 2 = t + t \bmod 2$,

$$\Delta^t O_{\tau,t}(|E|^{-(s+2)/2}) = \Delta^t O_{\tau,t}(|E|^{-(t+t \bmod 2)/2}) = O_{\tau,t}(|E|^{-(t+t \bmod 2)/2}).$$

Therefore,

$$\Delta^t \mathbb{P}_{E,\eta+k}(A) = O_{\tau,t,k}(|E|^{-(t+t \bmod)/2})$$

for all $t \geq 0$, and (51) follows. □

**4. High-order weighted sampling correlations.** For $E \subset \mathbf{N}$ and sets $(\mathbf{u} \cup \mathbf{v}) \subset E$ with $\mathbf{u} \cap \mathbf{v} = \varnothing$ and $\eta \geq |\mathbf{u}|$, we define the (conditional) measure $\mathbb{P}_{E,\eta}^{\mathbf{u},\mathbf{v}}$ supported on sets $\mathbf{d} \subset E$ of size $\eta$ with $\mathbf{d} \supset \mathbf{u}$ and $\mathbf{d} \cap \mathbf{v} = \varnothing$ by

$$\mathbb{P}_{E,\eta}^{\mathbf{u},\mathbf{v}}(\mathbf{D} = \mathbf{d}) = \mathbb{P}_{E,\eta}(\mathbf{D} = \mathbf{d} | \mathbf{D} \supset \mathbf{u}, \mathbf{D} \cap \mathbf{v} = \varnothing);$$

that is, $\mathbb{P}_{E,\eta}^{\mathbf{u},\mathbf{v}}$ is the measure $\mathbb{P}_{E,\eta}$ conditioned to contain every element of $\mathbf{u}$ but none of the elements of $\mathbf{v}$. The measures considered in the previous sections were the unconditioned special case

$$\mathbb{P}_{E,\eta} = \mathbb{P}_{E,\eta}^{\varnothing,\varnothing};$$

the $\varnothing, \varnothing$ superscript may be omitted. We define (commutative) differences $\Delta^B$ on the measure $\mathbb{P}_{E,\eta}^{\mathbf{u},\mathbf{v}}$ for $B \in E \setminus (\mathbf{u} \cup \mathbf{v})$ by

$$\Delta^B \mathbb{P}_{E,\eta}^{\mathbf{u},\mathbf{v}} = \mathbb{P}_{E,\eta}^{\mathbf{u} \cup B,\mathbf{v}} - \mathbb{P}_{E,\eta}^{\mathbf{u},\mathbf{v} \cup B}.$$

For $k \in \mathbf{N}$ the operators $\Delta^k$ will continue to be used in accordance with (30).

The following lemma gives some key properties of the conditional measure $\mathbb{P}_{E,\eta}^{\mathbf{u},\mathbf{v}}$, including a useful relation to its unconditional version $\mathbb{P}_{E,\eta}$.

LEMMA 4.1. *Let $\mathbf{u}, \mathbf{v}$ be disjoint subsets of $E$. For $H \subset E \setminus (\mathbf{u} \cup \mathbf{v})$,*

(52) $$\Delta^H \mathbb{P}_{E,\eta}^{\mathbf{u},\mathbf{v}} = \sum_{\alpha \cup \beta = H, \alpha \cap \beta = \varnothing} (-1)^{|\beta|} \mathbb{P}_{E,\eta}^{\mathbf{u} \cup \alpha, \mathbf{v} \cup \beta}.$$

*For $\mathbf{d} \subset E$ such that $\mathbf{d} \supset \mathbf{u}$ and $\mathbf{d} \cap \mathbf{v} = \varnothing$,*

(53) $$\mathbb{P}_{E,\eta}^{\mathbf{u},\mathbf{v}}(\mathbf{d}) = \mathbb{P}_{E \setminus (\mathbf{u} \cup \mathbf{v}), \eta - |\mathbf{u}|}(\mathbf{d} \setminus \mathbf{u}),$$

*and for $A \notin (\mathbf{u} \cup \mathbf{v})$ and $H \subset E \setminus (\mathbf{u} \cup \mathbf{v} \cup A)$,*

(54) $$\Delta^H \mathbb{P}_{E,\eta}^{\mathbf{u},\mathbf{v}}(A) = (-1)^{|H|} \Delta^{|H|} \mathbb{P}_{E \setminus H, \eta}^{\mathbf{u},\mathbf{v}}(A).$$

PROOF. Relation (52) can be shown by induction. By definition (4) and that of conditional probability, using $\mathbf{u} \subset \mathbf{d} \subset E \setminus \mathbf{v}$, we have

$$\mathbb{P}_{E,\eta}^{\mathbf{u},\mathbf{v}}(\mathbf{d}) = \frac{x_\mathbf{d}}{\sum_{\mathbf{w}: \mathbf{u} \subset \mathbf{w} \subset E, \mathbf{w} \cap \mathbf{v} = \varnothing, |\mathbf{w}| = \eta} x_\mathbf{w}}$$

$$= \frac{x_{\mathbf{d} \setminus \mathbf{u}}}{\sum_{\mathbf{w}: \mathbf{u} \subset \mathbf{w} \subset E, \mathbf{w} \cap \mathbf{v} = \varnothing, |\mathbf{w}| = \eta} x_{\mathbf{w} \setminus \mathbf{u}}},$$

since both $\mathbf{d}$ and $\mathbf{w}$ contain $\mathbf{u}$ and the factor $x_\mathbf{u}$, which appears in both $x_\mathbf{d} = x_{\mathbf{d} \setminus \mathbf{u}} x_\mathbf{u}$ and $x_\mathbf{w} = x_{\mathbf{w} \setminus \mathbf{u}} x_\mathbf{u}$, can be cancelled. Furthermore, because $\mathbf{u}, \mathbf{v}$ are disjoint, $\mathbf{w} \cap \mathbf{v} = \varnothing$ if and only if $(\mathbf{w} \setminus \mathbf{u}) \cap \mathbf{v} = \varnothing$, and when $\mathbf{u} \subset \mathbf{w}$ and $|\mathbf{w}| = \eta$, then $|\mathbf{w} \setminus \mathbf{u}| = \eta - |\mathbf{u}|$. Hence,

$$\{\mathbf{w} \setminus \mathbf{u} : \mathbf{u} \subset \mathbf{w} \subset E, \mathbf{w} \cap \mathbf{v} = \varnothing, |\mathbf{w}| = \eta\}$$
$$= \{\mathbf{w} \setminus \mathbf{u} : \mathbf{w} \setminus \mathbf{u} \subset E \setminus \mathbf{u}, (\mathbf{w} \setminus \mathbf{u}) \cap \mathbf{v} = \varnothing, |\mathbf{w} \setminus \mathbf{u}| = \eta - |\mathbf{u}|\}$$
$$= \{\mathbf{w} : \mathbf{w} \subset E \setminus \mathbf{u}, \mathbf{w} \cap \mathbf{v} = \varnothing, |\mathbf{w}| = \eta - |\mathbf{u}|\}$$
$$= \{\mathbf{w} : \mathbf{w} \subset E \setminus (\mathbf{u} \cup \mathbf{v}), |\mathbf{w}| = \eta - |\mathbf{u}|\}$$

and $\mathbb{P}^{\mathbf{u},\mathbf{v}}_{E,\eta}(\mathbf{d})$ equals

$$\frac{x_{\mathbf{d} \setminus \mathbf{u}}}{\sum_{\mathbf{w} \,:\, \mathbf{w} \subset E \setminus (\mathbf{u} \cup \mathbf{v}), |\mathbf{w}| = \eta - |\mathbf{u}|} x_\mathbf{w}} = \mathbb{P}_{E \setminus (\mathbf{u} \cup \mathbf{v}), \eta - |\mathbf{u}|}(\mathbf{d} \setminus \mathbf{u}).$$

This proves (53).

It suffices to prove (54) for $(\mathbf{u}, \mathbf{v}) = (\varnothing, \varnothing)$. First, note for $A \notin \alpha \cup \beta$,

$$\mathbb{P}^{\alpha,\beta}_{E,\eta}(A) = \sum_{\mathbf{d} \ni A, \alpha \subset \mathbf{d} \subset E, \mathbf{d} \cap \beta = \varnothing} \mathbb{P}^{\alpha,\beta}_{E,\eta}(\mathbf{d})$$
$$= \sum_{\mathbf{d} \ni A, \alpha \subset \mathbf{d} \subset E, \mathbf{d} \cap \beta = \varnothing} \mathbb{P}_{E \setminus (\alpha \cup \beta), \eta - |\alpha|}(\mathbf{d})$$
$$= \sum_{\mathbf{d} \ni A, \mathbf{d} \subset E \setminus (\alpha \cup \beta)} \mathbb{P}_{E \setminus (\alpha \cup \beta), \eta - |\alpha|}(\mathbf{d})$$

(55)
$$= \mathbb{P}_{E \setminus (\alpha \cup \beta), \eta - |\alpha|}(A);$$

hence, since $A \notin H$,

$$\Delta^H \mathbb{P}_{E,\eta}(A) = \sum_{\alpha \cup \beta = H, \alpha \cap \beta = \varnothing} (-1)^{|\beta|} \mathbb{P}^{\alpha,\beta}_{E,\eta}(A) \qquad \text{[by (52)]}$$
$$= \sum_{\alpha \cup \beta = H, \alpha \cap \beta = \varnothing} (-1)^{|\beta|} \mathbb{P}_{E \setminus (\alpha \cup \beta), \eta - |\alpha|}(A) \qquad \text{[by (55)]}$$
$$= \sum_{j=0}^{|H|} \sum_{\alpha \cup \beta = H, \alpha \cap \beta = \varnothing, |\alpha| = j} (-1)^{|H|-j} \mathbb{P}_{E \setminus H, \eta - j}(A)$$
$$= \sum_{j=0}^{|H|} \binom{|H|}{j} (-1)^{|H|-j} \mathbb{P}_{E \setminus H, \eta - j}(A)$$
$$= (-1)^{|H|} \sum_{j=0}^{|H|} \binom{|H|}{j} (-\Psi)^j \mathbb{P}_{E \setminus H, \eta}(A)$$
$$= (-1)^{|H|} (1 - \Psi)^{|H|} \mathbb{P}_{E \setminus H, \eta}(A)$$
$$= (-1)^{|H|} \Delta^{|H|} \mathbb{P}_{E \setminus H, \eta}(A). \qquad \square$$

In parallel to definition (31), for functions $f_E$ for which versions $f_E^{\alpha,\beta}$ are defined [such as $f_E = \mathbb{P}_{E,\eta}^{\mathbf{u},\mathbf{v}}(A)$], for $G \subset \mathbf{N}$ let

$$\Gamma_\mu^q(G) = \{f_E : \Delta^H f_E = O_{\mu,|H|}(|E|^{-(|H|+q+(|H|+q) \bmod 2)/2}) \text{ for all } H \cap G = \varnothing\}.$$

In parallel to Lemma 3.3, we have the following:

LEMMA 4.2. *Let $p \leq q$ be nonnegative integers, and suppose that $f_F \in \Gamma_\mu^p(P)$ and $g_F \in \Gamma_\mu^q(Q)$. Then*

$$P \subset Q \implies \Gamma_\mu^p(P) \supset \Gamma_\mu^q(Q),$$
$$af_F \in \Gamma_\mu^p(P), \ f_F + g_F \in \Gamma_\mu^p(P \cup Q),$$
$$P \cap H = \varnothing \implies \Delta^H f_F \in \Gamma_\mu^{p+|H|}(P \cup H),$$
$$f_{(F \setminus H) \cup G} \in \Gamma_\mu^p(P),$$
$$f_F g_F \in \Gamma_\mu^{p+q}(P \cup Q).$$

The proof, being parallel to that of Lemma 3.3, is omitted.

LEMMA 4.3. *Let Condition 3.1 hold and $\tau \in (0, 1/2]$. For $(E, \eta) \in \mathcal{E}_\tau$, $G \supset (\mathbf{u} \cup \mathbf{v})$ and $G \cap \{A\} = \varnothing$,*

(56) $$\mathbb{P}_{E,\eta}^{\mathbf{u},\mathbf{v}}(A) \in \Gamma_\tau^0(G \cup A),$$

(57) $$\mathbb{P}_{E,\eta}^{\mathbf{u},\mathbf{v}}(A)\mathbb{P}_{E,\eta}^{\mathbf{u},\mathbf{v}}(\bar{A}) \in \Gamma_\tau^0(G \cup A),$$

(58) $$\mathbb{P}_{E,\eta}^{\mathbf{u},\mathbf{v}}(A) - \mathbb{P}_{E,\eta}(A) \in \Gamma_\tau^1(G \cup A).$$

PROOF. For $H \subset E \setminus (G \cup A)$, by (54) and (53),

$$\Delta^H \mathbb{P}_{E,\eta}^{\mathbf{u},\mathbf{v}}(A) = (-1)^{|H|} \Delta^{|H|} \mathbb{P}_{E \setminus H, \eta}^{\mathbf{u},\mathbf{v}}(A)$$
$$= (-1)^{|H|} \Delta^{|H|} \mathbb{P}_{E \setminus (H \cup \mathbf{u} \cup \mathbf{v}), \eta - |\mathbf{u}|}(A).$$

The result (56) now follows by (51) of Corollary 3.1. Since $1 \in \Gamma_\tau^0(G \cup A)$, we have $\mathbb{P}_{E,\eta}^{\mathbf{u},\mathbf{v}}(\bar{A}) = 1 - \mathbb{P}_{E,\eta}^{\mathbf{u},\mathbf{v}}(A) \in \Gamma_\tau^0(G \cup A)$, and, hence, (57) using Lemma 4.2.

Next, if $B \in \mathbf{v} \neq \varnothing$,

$$\mathbb{P}_{E,\eta}^{\mathbf{u},\mathbf{v}}(A) - \mathbb{P}_{E,\eta}^{\mathbf{u} \cup B, \mathbf{v} \setminus B}(A) = -\Delta^B \mathbb{P}_{E,\eta}^{\mathbf{u},\mathbf{v} \setminus B}(A),$$

which is in $\Gamma_\tau^1(G \cup A)$ by (56). Iterating over all elements in $\mathbf{v}$ and using the fact that $\Gamma_\tau^1(G \cup A)$ is closed under addition, we obtain

(59) $$\mathbb{P}_{E,\eta}^{\mathbf{u},\mathbf{v}}(A) - \mathbb{P}_{E,\eta}^{\mathbf{u} \cup \mathbf{v}, \varnothing}(A) \in \Gamma_\tau^1(G \cup A).$$

Next, for $B \in \mathbf{u} \cup \mathbf{v}$,

$$\mathbb{P}_{E,\eta}(A) = \mathbb{P}_{E,\eta}(B)\mathbb{P}_{E,\eta}^{B,\varnothing}(A) + \mathbb{P}_{E,\eta}(\bar{B})\mathbb{P}_{E,\eta}^{\varnothing,B}(A)$$

$$= (1 - \mathbb{P}_{E,\eta}(\bar{B}))\mathbb{P}_{E,\eta}^{B,\varnothing}(A) + \mathbb{P}_{E,\eta}(\bar{B})\mathbb{P}_{E,\eta}^{\varnothing,B}(A)$$

$$= \mathbb{P}_{E,\eta}^{B,\varnothing}(A) - \mathbb{P}_{E,\eta}(\bar{B})(\mathbb{P}_{E,\eta}^{B,\varnothing}(A) - \mathbb{P}_{E,\eta}^{\varnothing,B}(A)).$$

Rearranging,

$$\mathbb{P}_{E,\eta}^{B,\varnothing}(A) - \mathbb{P}_{E,\eta}(A) = \mathbb{P}_{E,\eta}(\bar{B})\Delta^B \mathbb{P}_{E,\eta}(A).$$

Since $\mathbb{P}_{E,\eta}(\bar{B}) \in \Gamma_\tau^0(G \cup A)$ and $\Delta^B \mathbb{P}_{E,\eta}(A) \in \Gamma_\tau^1(G \cup A)$, their product is in $\Delta^B \mathbb{P}_{E,\eta}(A) \in \Gamma_\tau^1(G \cup A)$ by (4.2) and, therefore, $\mathbb{P}_{E,\eta}^{B,\varnothing}(A) - \mathbb{P}_{E,\eta}(A) \in \Gamma_\tau^1(G \cup A)$. Iterating over all $B \in \mathbf{u} \cup \mathbf{v}$ and using the fact that $\Gamma_\tau^1(G \cup A)$ is closed under addition, we have

$$\mathbb{P}_{E,\eta}^{\mathbf{u} \cup \mathbf{v}, \varnothing}(A) - \mathbb{P}_{E,\eta}(A) \in \Gamma_\tau^1(G \cup A),$$

and now by (59) and the closure property of $\Gamma_\tau^1(G \cup A)$ (58) follows. □

For short, write

$$p_A^{\mathbf{u},\mathbf{v}} = \mathbb{P}_{E,\eta}^{\mathbf{u},\mathbf{v}}(A) \quad \text{and} \quad q_A^{\mathbf{u},\mathbf{v}} = 1 - p_A^{\mathbf{u},\mathbf{v}};$$

as usual, for $(\mathbf{u}, \mathbf{v}) = (\varnothing, \varnothing)$, we omit the superscripts.

LEMMA 4.4. *For any random variable $V$ and $A \notin \mathbf{u} \cup \mathbf{v}$,*

$$\mathbb{E}_{E,\eta}^{\mathbf{u},\mathbf{v}}((I_A - p_A)V) = ((p_A^{\mathbf{u},\mathbf{v}} - p_A) + p_A^{\mathbf{u},\mathbf{v}} q_A^{\mathbf{u},\mathbf{v}} \Delta^A)\mathbb{E}_{E,\eta}^{\mathbf{u},\mathbf{v}}(V).$$

PROOF. Adding and subtracting $p_A^{\mathbf{u},\mathbf{v}}$, we have

$$\mathbb{E}_{E,\eta}^{\mathbf{u},\mathbf{v}}((I_A - p_A)V) = \mathbb{E}_{E,\eta}^{\mathbf{u},\mathbf{v}}((I_A - p_A^{\mathbf{u},\mathbf{v}})V) + (p_A^{\mathbf{u},\mathbf{v}} - p_A)\mathbb{E}_{E,\eta}^{\mathbf{u},\mathbf{v}}(V)$$

and

$$\mathbb{E}_{E,\eta}^{\mathbf{u},\mathbf{v}}((I_A - p_A^{\mathbf{u},\mathbf{v}})V)$$
$$= p_A^{\mathbf{u},\mathbf{v}} \mathbb{E}_{E,\eta}^{\mathbf{u} \cup A, \mathbf{v}}((1 - p_A^{\mathbf{u},\mathbf{v}})V) + (1 - p_A^{\mathbf{u},\mathbf{v}})\mathbb{E}_{E,\eta}^{\mathbf{u}, \mathbf{v} \cup A}(-p_A^{\mathbf{u},\mathbf{v}} V)$$
$$= p_A^{\mathbf{u},\mathbf{v}} q_A^{\mathbf{u},\mathbf{v}} (\mathbb{E}_{E,\eta}^{\mathbf{u} \cup A, \mathbf{v}}(V) - \mathbb{E}_{E,\eta}^{\mathbf{u}, \mathbf{v} \cup A}(V))$$
$$= p_A^{\mathbf{u},\mathbf{v}} q_A^{\mathbf{u},\mathbf{v}} \Delta^A \mathbb{E}_{E,\eta}^{\mathbf{u},\mathbf{v}}(V). \qquad \square$$

THEOREM 4.1. *Let Condition 3.1 hold and $(E, \eta) \in \mathcal{E}_\tau$ for $\tau \in (0, 1/2]$. If $\mathbf{u}, \mathbf{v}$ are subsets of $E$ with $G \supset (\mathbf{u} \cup \mathbf{v})$ and $V$ is a random variable such that*

$$\mathbb{E}_{E,\eta}^{\mathbf{u},\mathbf{v}}(V) \in \Gamma_\tau^q(G),$$

*then for $G \cap H = \varnothing$,*

$$\mathbb{E}_{E,\eta}^{\mathbf{u},\mathbf{v}}\left(\prod_{A\in H}(I_A - p_A)V\right) \in \Gamma_\tau^{q+|H|}(G \cup H).$$

*In particular, when $V = 1$ and $G = u = v = \varnothing$, since $1 \in \Gamma_\tau^0(\varnothing)$, we have*

$$\mathrm{Corr}(H) \equiv \mathbb{E}_{E,\eta}\left(\prod_{A\in H}(I_A - p_A)\right) \in \Gamma_\tau^{|H|}(H),$$

*and, therefore, in particular,*

$$\mathrm{Corr}(H) = O_{\tau,|H|}(|E|^{-(|H|+|H|\bmod 2)/2}).$$

PROOF. For $H = \{A\}$, by Lemma 4.2,

$$\Delta^A \mathbb{E}_{E,\eta}^{\mathbf{u},\mathbf{v}}(V) \in \Gamma_\tau^{q+1}(G \cup A);$$

since $p_A^{\mathbf{u},\mathbf{v}} q_A^{\mathbf{u},\mathbf{v}} \in \Gamma_\tau^0(G)$, using Lemma 4.2 again yields

$$p_A^{\mathbf{u},\mathbf{v}} q_A^{\mathbf{u},\mathbf{v}} \Delta^A \mathbb{E}_{E,\eta}^{\mathbf{u},\mathbf{v}}(V) \in \Gamma_\tau^{q+1}(G \cup A).$$

Since

$$p_A^{\mathbf{u},\mathbf{v}} - p_A \in \Gamma_\tau^1(G \cup A),$$

we also have that

$$(p_A^{\mathbf{u},\mathbf{v}} - p_A)\mathbb{E}_{E,\eta}^{\mathbf{u},\mathbf{v}}(V) \in \Gamma_\tau^{q+1}(G \cup A).$$

The result for $H = \{A\}$ now follows from Lemma 4.4, and then, in general, by induction. □

We close this section with some results which will be useful in Section 5.

COROLLARY 4.1. *Under the hypotheses for Theorem 4.1, for $A, B, C$ distinct,*

$$p_A = p_{A,\boldsymbol{\lambda}} + O_\tau(|E|^{-1}),$$
$$\mathbb{E}_{E,\eta}(I_A - p_A)(I_B - p_B) = -p_{A,\boldsymbol{\lambda}} q_{A,\boldsymbol{\lambda}} p_{B,\boldsymbol{\lambda}} q_{B,\boldsymbol{\lambda}} v_{E,\boldsymbol{\lambda}}^{-2} + o_\tau(|E|^{-1}),$$
$$\mathbb{E}_{E,\eta}(I_A - p_A)^2(I_B - p_B) = O_\tau(|E|^{-1}),$$
$$\mathbb{E}_{E,\eta}(I_A - p_A)(I_B - p_B)(I_C - p_C) = O_\tau(|E|^{-2}).$$

PROOF. The first claim is a consequence of (50). With $F = E \setminus (A \cup B)$,

$$
\begin{aligned}
\mathbb{E}_{E,\eta}(I_A - p_A)&(I_B - p_B) \\
&= p_A q_A \Delta^A \mathbb{E}_{E,\eta}(I_B - p_B) && \text{(by Lemma 4.4)} \\
&= p_A q_A \Delta^A \mathbb{P}_{E,\eta}(B) && [\text{since } \mathbb{E}_{E,\eta}^{\mathbf{u},\mathbf{v}}(p_B) = p_B \text{ for all } \mathbf{u},\mathbf{v}] \\
&= -p_A q_A \Delta \mathbb{P}_{E \setminus A,\eta}(B) && [\text{by (54) of Lemma 4.1}] \\
&= -p_A q_A p_{B,\boldsymbol{\lambda}} q_{B,\boldsymbol{\lambda}} \mathcal{I}_{F,\boldsymbol{\lambda},2}^0 + O_\tau(|E|^{-2}) && [\text{by (50) of Theorem 3.1}] \\
&= -p_{A,\boldsymbol{\lambda}} q_{A,\boldsymbol{\lambda}} p_{B,\boldsymbol{\lambda}} q_{B,\boldsymbol{\lambda}} \mathcal{I}_{F,\boldsymbol{\lambda},2}^0 + O_\tau(|E|^{-2}) && [\text{by (19) and (50)}] \\
&= -p_{A,\boldsymbol{\lambda}} q_{A,\boldsymbol{\lambda}} p_{B,\boldsymbol{\lambda}} q_{B,\boldsymbol{\lambda}} v_{F,\boldsymbol{\lambda}}^{-2} + o_\tau(|E|^{-1}) && [\text{by (19)}] \\
&= -p_{A,\boldsymbol{\lambda}} q_{A,\boldsymbol{\lambda}} p_{B,\boldsymbol{\lambda}} q_{B,\boldsymbol{\lambda}} v_{E,\boldsymbol{\lambda}}^{-2} + o_\tau(|E|^{-1}) && \text{(by Lemma 3.1).}
\end{aligned}
$$

Further,

$$
\begin{aligned}
\mathbb{E}_{E,\eta}(I_A - p_A)^2 &(I_B - p_B) \\
&= p_A(1-p_A)^2 \mathbb{E}_{E,\eta}^{A,\varnothing}(I_B - p_B) + (1-p_A)p_A^2 \mathbb{E}_{E,\eta}^{\varnothing,A}(I_B - p_B) \\
&= p_A(1-p_A)((1-p_A)(p_B^{A,\varnothing} - p_B) + p_A(p_B^{\varnothing,A} - p_B)) \\
&= O_\tau(|E|^{-1}) \qquad \text{(by Lemma 4.3),}
\end{aligned}
$$

and the final claim is immediate from Theorem 4.1. □

**5. Application: asymptotics for conditional and unconditional logistic odds ratio estimators.** In this section the theory developed in the previous sections is used to provide an asymptotic theory for the maximum likelihood conditional and unconditional logistic regression odds ratio estimators, $\hat{\boldsymbol{\beta}}_N$ and $\tilde{\boldsymbol{\beta}}_N$, under the nested case-control model. Conditions 5.1 and 5.2 ensure the asymptotic stability and nondegeneracy of data in the study base, which is sampled using schemes satisfying Condition 5.3. Lemma 5.1 shows how stability in the study base leads to stability in probability for case-control samples $E$. Theorems 5.1 and 5.2 give the consistency and asymptotic normality of $\hat{\boldsymbol{\beta}}_N$ and $\tilde{\boldsymbol{\beta}}_N$. We first consider asymptotically stable covariates in $\mathcal{R}$ and then specialize to the i.i.d. case. Previously, the weights $x_A, A \in E$ were considered fixed, but here even if $x_j, j \in \mathcal{R}$ are fixed, the values $x_A, A \in E$ arrive in $E$ through random failure and control sampling. Suppressing explicit dependence on $\boldsymbol{\beta}_0$ and (as usual) on $\mathcal{R}$ and its size $N$, we indicate the study base model $\mathbb{P}_{\lambda_0,\beta_0}$ given in (1) by $\mathbb{P}$, and continue to denote the conditional distributions given $E, \eta$ by $\mathbb{P}_{E,\eta}$.

The first two conditions are on the stability of the study base data.

CONDITION 5.1. For all $\delta \in (0,1)$ there exists $C$ such that for all $N \geq 1$,

$$\frac{1}{N} \sum_{j \in \mathcal{R}} \mathbf{1}(|\mathbf{z}_j| \leq C) \geq 1 - \delta \tag{60}$$

and with $p_j$ given by (1) with $x = x_j$,

$$\frac{1}{N} \sum_{j \in \mathcal{R}} p_j \to p \quad \text{as } N \to \infty.$$

Clearly, we then have $\eta/N \xrightarrow{p} p$ as $N \to \infty$ by independence of the failure indicators. Furthermore, $p \in (0,1)$, since with $C$ corresponding to any $\delta \in (0,1)$ in (60),

$$\frac{1}{N} \sum_{j \in \mathcal{R}} p_j \geq \left( \inf_{|\mathbf{z}| \leq C} p(\mathbf{z}) \right) \frac{1}{N} \sum_{j \in \mathcal{R}} \mathbf{1}(|\mathbf{z}_j| \leq C) \geq \left( \inf_{|\mathbf{z}| \leq C} p(\mathbf{z}) \right)(1 - \delta),$$

which by Condition 1.1 is strictly positive for all $N \geq 1$; likewise for $q_j$, where $p_j + q_j = 1$. For $u_j, j \in \mathcal{R}$, let

$$\overline{u}_N = \frac{1}{N} \sum_{j \in \mathcal{R}} u_j \quad \text{and} \quad \overline{\overline{u}} = \sup_{N \geq 1} \overline{u}_N. \tag{61}$$

We say $u_j$ is asymptotically stable in mean if $\overline{u}_N \to \overline{u}$ for some $\overline{u}$ as $N \to \infty$, asymptotically dominated in mean if $\overline{\overline{|u|}} < \infty$, and $u_j(\boldsymbol{\beta})$ uniformly asymptotically dominated in mean if there exists a neighborhood $\mathcal{B}_0 \subset \mathcal{B}$ containing $\boldsymbol{\beta}_0$, and $v_j$ asymptotically dominated in mean, such that $|u_j(\boldsymbol{\beta})| \leq v_j$ for all $\boldsymbol{\beta} \in \mathcal{B}_0$. For a continuous function $w: [0, \infty) \to \mathbb{R}$ with $\lim_{x \to \infty} w(x) = L \in (-\infty, \infty)$, we say $u_j$ is $w$-stable if $|u_j|^2$ is asymptotically dominated in mean and for all $\lambda \in [0, \infty]$, $u_j w(\lambda x_j) p_j$ and $u_j w(\lambda x_j) q_j$ are asymptotically stable in mean. In what follows, we omit the specification "in mean."

CONDITION 5.2. 1 is $x/(1+x)$ stable, $\mathbf{z}_j^{\otimes k}$ is $x/(1+x)^2$ stable for $k = 0, 1, 2$, $|\mathbf{z}_j|^3, |\mathbf{z}_j'|^2$ and $|\mathbf{z}_j''|^2$ are uniformly asymptotically dominated, and

$$\liminf_{N \to \infty} \inf_{|\mathbf{a}|=1} \mathbf{a}^\mathsf{T} \left( \frac{1}{N} \sum_{j \in \mathcal{R}} (\mathbf{y}_j - \bar{\mathbf{y}}_N)^{\otimes 2} \right) \mathbf{a} > 0, \tag{62}$$

for $\mathbf{y}_j^\mathsf{T} = (1, \mathbf{z}_j^\mathsf{T})$.

The next condition is on the sampling design.

CONDITION 5.3. For some $f \in (0,1)$,

$$\frac{\eta}{|E|} \xrightarrow{p} f \quad \text{as } N \to \infty, \tag{63}$$

for $B_j = \mathbf{1}(j \in E)$, uniformly over $j \in \mathcal{R}$,

$$\mathbb{E}(B_j | j \notin \mathbf{D}) \to \rho_f = \frac{(1-f)p}{(1-p)f} \tag{64}$$

and uniformly over all $j \neq k$ in $\mathcal{R}$,

$$\mathrm{Cov}((B_j, B_k) | j \notin \mathbf{D}, k \notin \mathbf{D}) \to 0 \qquad \text{as } N \to \infty.$$

For $f$ as in Condition 5.3, set $\tau = (1/2)\min(f, 1-f)$ for application of Corollary 4.1. The connection between properties of $u_j$ on $\mathcal{R}$ and their corresponding in probability versions on $E$ is made explicit by Lemma 5.1. We say $g_E(\lambda)$ converges uniformly in probability to $g(\lambda)$ if $\sup_{\lambda \in [0,\infty)} |g_E(\lambda) - g(\lambda)| \xrightarrow{p} 0$ as $N \to \infty$.

LEMMA 5.1. *Assume Conditions* 1.1, 5.1 *and* 5.3 *hold.*

(a) *For all $\delta \in (0,1)$, there exists $\varepsilon \in (0,1)$ such that for all $N \geq 1$,*

$$\mathbb{P}\left(\frac{1}{|E|} \sum_{A \in E} \mathbf{1}(x_A \in [\varepsilon, \varepsilon^{-1}]) \geq 1 - \delta\right) \geq 1 - \delta. \tag{65}$$

*If $u_j$ is asymptotically dominated, then for all $\delta \in (0,1)$, there exists $K$ such that*

$$\mathbb{P}\left(\frac{1}{|E|} \sum_{A \in E} |U_A| \leq K\right) \geq 1 - \delta \qquad \text{for all } N \geq 1. \tag{66}$$

(b) *If $|u_j|^2$ is asymptotically dominated, then $\mathrm{Var}(N^{-1} \sum_{A \in E} U_A) \to 0$. If, in addition, $u_j p_j$ and $u_j q_j$ are asymptotically stable, then*

$$\frac{1}{|E|} \sum_{A \in E} U_A \xrightarrow{p} f\frac{\overline{up}}{p} + (1-f)\frac{\overline{uq}}{q}$$

$$= \frac{1-f}{1-p}\left[u\frac{1 + \rho_f^{-1}\lambda_0 x}{1 + \lambda_0 x}\right] \qquad \text{as } N \to \infty, \tag{67}$$

*with $\rho_f$ given in* (64).

(c) *If $u_j$ is w-stable, then*

$$g_E^U(\lambda) = \frac{1}{|E|} \sum_{A \in E} U_A w(\lambda x_A) \tag{68}$$

*converges in probability uniformly to a continuous limit $g^U(\lambda)$ as $N \to \infty$, having form* (67) *with $u$ replaced by $uw(\lambda x)$. Hence, additionally, under Condition* 5.2,

$$h_E(\lambda) = \frac{1}{|E|} \sum_{A \in E} p_{A,\lambda} \quad \text{and} \quad e_{k,E}(\lambda) = \frac{p}{f}\frac{1}{|E|} \sum_{A \in E} \mathbf{z}_A^{\otimes k} p_{A,\lambda} q_{A,\lambda}$$

*converge uniformly in probability to continuous functions* $h(\lambda)$ *and* $e_k(\lambda)$ *for* $k = 0, 1, 2$ *with form* (67).

(d) *The limit function* $h(\lambda)$ *in part* (c) *strictly increases from 0 to 1 as* $\lambda$ *increases from 0 to* $\infty$. *For* $f \in (0,1)$,

$$\lambda_f = \rho_f^{-1} \lambda_0$$

*is the unique solution to* $h(\lambda_f) = f$, *and*

(69) $$e_k \equiv e_k(\lambda_f) = \overline{[\mathbf{z}^{\otimes k} q_{\lambda_f} p_{\lambda_0}]}.$$

*With* $h_E(\boldsymbol{\lambda}) = \eta/|E|$, *we have*

$$\boldsymbol{\lambda} \xrightarrow{p} \lambda_f.$$

(e) *if* $|u_j|^2$ *and* $|v_j|^2$ *are (uniformly) asymptotically dominated, then*

$$\frac{1}{|E|} \operatorname{Cov}_{E,\eta}\left(\sum_{A \in \mathbf{D}} U_A, \sum_{A \in \mathbf{D}} V_A\right)$$

$$= \frac{1}{|E|}\left(\sum_{A \in E} U_A V_A p_A q_A + \sum_{A \neq B} U_A V_B (p_{AB} - p_A p_B)\right)$$

*is (uniformly) bounded in probability. If* $1, u_j, v_j$ *and* $u_j v_j$ *are* $w(x) = x/(1+x)^2$*-stable, then with* $g^U$ *the limit of* $g_E^U$ *given in* (68),

$$|E|^{-1} \operatorname{Cov}_{E,\eta}\left(\sum_{A \in \mathbf{D}} U_A, \sum_{A \in \mathbf{D}} V_A\right)$$

$$\xrightarrow{p} g^{U \otimes V}(\lambda_f) - pf^{-1} g^U(\lambda_f) \otimes g^V(\lambda_f)/e_0(\lambda_f).$$

(f) *If* $|u_j|^2$ *is (uniformly) asymptotically dominated, then (uniformly)*

$$\frac{1}{N} \sum_{A \in E} U_A(I_A - p_A) \xrightarrow{p} 0.$$

(g) *If nonnegative weights* $w_j$ *are asymptotically dominated, for all* $\delta \in (0,1)$, *there exists* $\varepsilon > 0$ *such that*

$$\inf_{N \geq 1} \frac{1}{N} \sum_{j \in \mathcal{R}} \mathbf{1}(w_j \geq \varepsilon) \geq 1 - \delta,$$

$|u_j|^3$ *is asymptotically dominated and*

(70) $$\liminf_{N \to \infty} \inf_{|\mathbf{a}|=1} \mathbf{a}^\mathsf{T} \Gamma_N \mathbf{a} > 0 \qquad \textit{where } \Gamma_N = \frac{1}{N} \sum_{j \in \mathcal{R}} (u_j - \overline{u}_N)^{\otimes 2},$$

*then this same lower bound holds for*

$$\Gamma_{N,w} = \sum_{j \in \mathcal{R}} (u_j - \overline{u}_{N,w})^{\otimes 2} \frac{w_j}{\sum_{k \in \mathcal{R}} w_k}, \qquad \text{where } \overline{u}_{N,w} = \sum_{j \in \mathcal{R}} u_j \frac{w_j}{\sum_{k \in \mathcal{R}} w_k}.$$

*In particular, under Condition* 5.2 *with* $e_k, k = 0, 1, 2$ *given in* (69),

(71) $$\Sigma = e_2 - e_0^{-1} e_1^{\otimes 2} \text{ is positive definite.}$$

PROOF. By considering coordinates, we will assume when convenient that $u_j \in \mathbb{R}$. To show (65), note that

$$\frac{1}{|E|} \sum_{A \in E} \mathbf{1}(x_A \notin [\varepsilon, \varepsilon^{-1}]) \leq \frac{N}{|E|} \frac{1}{N} \sum_{j \in \mathcal{R}} \mathbf{1}(x_A \notin [\varepsilon, \varepsilon^{-1}]).$$

By Conditions 1.1 and 5.1, $N^{-1} \sum_{j \in \mathcal{R}} \mathbf{1}(x_j \notin [\varepsilon, \varepsilon^{-1}])$ can be made arbitrarily small for all $N \geq 1$ by choosing $\varepsilon \in (0,1)$ sufficiently small. Now by (63) of Condition 5.3 and $\eta/N \xrightarrow{p} p$, we have $|E|/N \xrightarrow{p} p/f$, and (65) follows. Claim (66) follows from

$$\frac{1}{N} \sum_{A \in E} |U_A| \leq \frac{1}{N} \sum_{j \in \mathcal{R}} |u_j| \leq \overline{\overline{u}},$$

Chebyshev's inequality, and $|E|/N \xrightarrow{p} p/f$.

For (b) note that $E$ is comprised of the set of failures $\mathbf{D}$ from $\mathcal{R}$ and a sample from the complement $\mathcal{R} \setminus \mathbf{D}$. Hence,

(72) $$\frac{1}{N} \sum_{A \in E} U_A = \frac{1}{N} \sum_{j \in \mathcal{R}} u_j \mathbf{1}(j \in \mathbf{D}) + \frac{1}{N} \sum_{j \in \mathcal{R}} u_j \mathbf{1}(j \in E \setminus \mathbf{D}).$$

Apply $\operatorname{Var}(X+Y) \leq 2(\operatorname{Var}(X) + \operatorname{Var}(Y))$ on the right-hand side of (72). For the first term, by independence, $N^{-1} \operatorname{Var}(\sum_{j \in \mathcal{R}} u_j \mathbf{1}(j \in \mathbf{D})) \leq \overline{\overline{u^2}}$.

The indicators in the second term of (72) may not be independent. Write its variance as the sum of the diagonal term

$$\frac{1}{N^2} \sum_{j \in \mathcal{R}} u_j^2 q_j \mathbb{E}(B_j | j \notin \mathbf{D})(1 - q_j \mathbb{E}(B_j | j \notin \mathbf{D})) \leq N^{-1} \overline{\overline{u^2}} \to 0,$$

and the covariance term, with $c_N = \max_{j,k} |\operatorname{Cov}(B_j, B_k | j \notin \mathbf{D}, k \notin \mathbf{D})|$,

$$\frac{1}{N^2} \left| \sum_{j \neq k} u_j u_k q_j q_k \operatorname{Cov}(B_j, B_k | j \notin \mathbf{D}, k \notin \mathbf{D}) \right| \leq c_N (\overline{\overline{u}})^2 \to 0;$$

hence, $\operatorname{Var}(N^{-1} \sum_{A \in E} U_A) \to 0$.

From (72),

(73) $$\mathbb{E}\left( \frac{1}{N} \sum_{A \in E} U_A \right) = \frac{1}{N} \sum_{j \in \mathcal{R}} u_j p_j + \frac{1}{N} \sum_{j \in \mathcal{R}} u_j q_j \mathbb{E}(B_j | j \notin \mathbf{D}).$$

Using (64) of Condition 5.3 and the fact that $u_j$ is dominated, the limit of the difference between the expectation (73) and $\overline{up}_N + \overline{uq}_N \rho_f$ is zero. The first equality in (67) now follows from the first part, the stability conditions on $u_j$, that $N/|E| \xrightarrow{p} f/p$ and the definition of $\rho_f$. The second equality now follows from (1), which gives the identity

$$f\frac{\overline{up}_N}{p} + (1-f)\frac{\overline{uq}_N}{q} = \frac{1-f}{1-p}\frac{1}{N}\sum_{j\in\mathcal{R}} u_j(\rho_f^{-1} p_j + q_j).$$

Turning to (c), since $w$ is continuous with finite limit at infinity, and since the stability conditions hold in $[0,\infty]$, without loss of generality, through the mapping $\lambda \to \lambda/(1+\lambda)$ say, it suffices to consider $\lambda \in [0,1]$. Let $u(j,\lambda)$ stand for either $u_j w(\lambda x_j) p_j$ or $u_j w(\lambda x_j) q_j$. Since $\|w\| = \sup_{\lambda \in [0,1]} |w(\lambda)| < \infty$, we have $|u(j,\lambda)|^2 \leq \|w\|^2 |u_j|^2$ and part (b) now shows that for all $\lambda$, $g_E(\lambda)$ converges in probability to $g(\lambda)$ having form claimed. It remains to show that the limit is continuous and that the convergence is uniform.

Let $\delta \in (0,1)$ be given. Since

$$\overline{U_E^2} \equiv \frac{1}{|E|}\sum_{A\in E} U_A^2 \leq \frac{N}{|E|}\overline{\overline{u^2}} \xrightarrow{p} \frac{f}{p}\overline{\overline{u^2}},$$

there is $M \geq 1$ such that for all $E$,

(74) $\qquad P(|\overline{U_E^2}| \leq K) \geq 1 - \delta/6, \qquad \text{where } K = Mf\overline{\overline{u^2}}/p.$

Assume for nontriviality that $\|w\|$ and $\overline{\overline{u^2}}$ are positive. Setting for short

$$\mathbf{1}_A(\varepsilon) = \mathbf{1}(x_A \notin [\varepsilon, \varepsilon^{-1}]),$$

and using notation as in (61), by part (a), there exists $\varepsilon \in (0,1)$ such that for all $E$,

$$\mathbb{P}(\overline{\mathbf{1}}_E(\varepsilon) \leq \delta^2/(16\|w\|^2 K)) \geq 1 - \delta/6.$$

Writing $g_E(\lambda)$ for $g_E^U(\lambda)$, let

(75) $\qquad g_E(\lambda) = g_E^{\triangleright\triangleleft}(\lambda) + g_E^{\triangleleft\triangleright}(\lambda),$

where

$$g_E^{\triangleright\triangleleft}(\lambda) = \frac{1}{|E|}\sum_{A\,:\,x_A\in[\varepsilon,\varepsilon^{-1}]} U_A w(\lambda x_A)$$

and

$$g_E^{\triangleleft\triangleright}(\lambda) = \frac{1}{|E|}\sum_{A\,:\,x_A\notin[\varepsilon,\varepsilon^{-1}]} U_A w(\lambda x_A).$$

Now applying the Cauchy–Schwarz inequality, with probability at least $1 - \delta/3$,

$$\sup_{0 \leq \lambda_2, \lambda_1 \leq 1} |g_E^{\triangleleft \triangleright}(\lambda_2) - g_E^{\triangleleft \triangleright}(\lambda_1)|$$

$$\leq 2 \sup_{0 \leq \lambda \leq 1} |g_E^{\triangleleft \triangleright}(\lambda)|$$

$$\leq 2\|w\| \frac{1}{|E|} \sum_{A \in E} |U_A| \mathbf{1}_A(\varepsilon) \leq 2\|w\| (\overline{U_E^2} \, \overline{\mathbf{1}}_E(\varepsilon))^{1/2} \leq \frac{\delta}{2}.$$

Since $w$ is uniformly continuous on $[0,1]$, there exists $\tau > 0$ such that

$$\text{if} \quad |y - x| < \tau/\varepsilon \quad \text{then } |w(y) - w(x)| < \delta/(2K^{1/2}).$$

In particular,

when $x_j \leq \varepsilon^{-1}$, if $|\lambda_2 - \lambda_1| < \tau$ then $|w(\lambda_2 x_j) - w(\lambda_1 x_j)| < \delta/(2K^{1/2})$.

Hence, by (74), with probability at least $1 - \delta/6$,

$$|g_E^{\triangleright \triangleleft}(\lambda_2) - g_E^{\triangleright \triangleleft}(\lambda_1)| \leq \frac{1}{|E|} \sum_{A \,:\, x_A \in [\varepsilon, \varepsilon^{-1}]} |U_A| |w(\lambda_2 x_A) - w(\lambda_1 x_A)|$$

$$\leq \frac{\delta \sqrt{|\overline{U_E^2}|}}{2K^{1/2}} \leq \frac{\delta}{2}.$$

Now by (75), for every $\delta$ there is a $\tau$ such that for all $E$,

$$\mathbb{P}\left( \sup_{|\lambda_2 - \lambda_1| \leq \tau} |g_E(\lambda_2) - g_E(\lambda_1)| \leq \delta \right) \geq 1 - \delta/2,$$

and taking limits, $\sup_{|\lambda_2 - \lambda_1| \leq \tau} |g(\lambda_2) - g(\lambda_1)| \leq \delta$; hence, $g(\lambda)$ is continuous. Letting $F_1, \ldots, F_M$ be a finite subcover of $[0,1]$ taken from the open cover of all open sub-intervals of length $2\tau$ and setting $\lambda_j$ to be the center of the interval $F_j$, there exists $N_0$ such that for $|E| \geq N_0$,

$$\mathbb{P}(|g_E(\lambda_j) - g(\lambda_j)| \leq \delta) \geq 1 - \delta/2, \qquad j = 1, \ldots, M.$$

Now (c) is finished, since for any $\lambda$, there exists $\lambda_j$ with $|\lambda - \lambda_j| < \tau$, and

$$|g_E(\lambda) - g(\lambda)| \leq |g_E(\lambda) - g_E(\lambda_j)| + |g_E(\lambda_j) - g(\lambda_j)| + |g(\lambda_j) - g(\lambda)|,$$

and so for all $\delta$ there exists $N_0$ such that

$$\text{for all } |E| \geq N_0 \qquad \mathbb{P}\left( \sup_{\lambda \in [0,1]} |g_E(\lambda) - g(\lambda)| \leq 3\delta \right) \geq 1 - \delta.$$

To show (d), as in part (a), for given $\delta \in (0,1)$, there exists $\varepsilon \in (0,1)$ such that for all $N \geq 1$,

$$\frac{1}{N} \sum_{j \in \mathcal{R}} \mathbf{1}(x_j \in [\varepsilon, \varepsilon^{-1}]) \geq 1 - \delta.$$

Let $0 \le \lambda_1 < \lambda_2 < \infty$ and set

$$\gamma = \inf_{x_j \in [\varepsilon, \varepsilon^{-1}]} (p_{j,\lambda_2} - p_{j,\lambda_1}) p_j,$$

which is strictly positive. Since $p_{j,\lambda}$ is nondecreasing in $\lambda$, for all $N \ge 1$,

$$\overline{p_{\lambda_2} p}_N - \overline{p_{\lambda_1} p}_N \ge \frac{1}{N} \sum_{j \,:\, x_j \in [\varepsilon, \varepsilon^{-1}]} (p_{j,\lambda_2} - p_{j,\lambda_1}) p_j \ge \gamma(1-\delta),$$

and, hence, $\overline{p_{\lambda_2} p} > \overline{p_{\lambda_1} p}$; similarly, $\overline{p_{\lambda_2} q} > \overline{p_{\lambda_1} q}$. As the form of the limit function $h$ is given by (67) with $u_j = p_{j,\lambda}$, $h$ strictly increases from 0 to 1 as $\lambda$ increases from 0 to $\infty$. By continuity, for every $f \in (0,1)$ there exists a unique $\lambda_f$ such that $h(\lambda_f) = f$.

Next, note that setting $\lambda_f = \rho_f^{-1} \lambda_0$, we have

$$p_{j,\lambda_f}\left( \frac{1 + \rho_f^{-1} \lambda_0 x_j}{1 + \lambda_0 x_j} \right) = \rho_f^{-1} p_{j,\lambda_0},$$

which by (67) gives

$$h_E(\lambda_f) \xrightarrow{p} \frac{1-f}{1-p} \rho_f^{-1} p = f$$

and the claimed representation of $e_k(\lambda_f)$.

Last, since $h_E(\lambda(E,\eta)) = \eta/|E|$,

$$h(\boldsymbol{\lambda}) - h(\lambda_f) = h(\lambda(E,\eta)) - h_E(\lambda(E,\eta)) - \left( f - \frac{\eta}{|E|} \right) \xrightarrow{p} 0,$$

we have

$$\boldsymbol{\lambda} \xrightarrow{p} \lambda_f \qquad \text{as } |E| \to \infty,$$

since $h(\lambda)$ is continuous and strictly increasing.

For (e), by Corollary 4.1, the correspondence

$$v_{E,\lambda}^{-2} = \frac{1}{|E|} \frac{p}{f} e_{0,E}^{-1}(\lambda),$$

and the (uniform) domination assumed on $u_j, v_j$, we have that

$$\frac{1}{|E|} \sum_{A \in E} U_A V_A p_A q_A + \frac{1}{|E|} \sum_{A \ne B} U_A V_B (p_{AB} - p_A p_B)$$

is in probability (uniformly) within $o_\tau(1)$ of

$$\frac{1}{|E|} \sum_{A \in E} U_A V_A p_{A,\boldsymbol{\lambda}} q_{A,\boldsymbol{\lambda}} - \frac{1}{|E|^2} \frac{p}{f} \sum_{A,B \in E, A \ne B} U_A V_B p_{A,\boldsymbol{\lambda}} q_{A,\boldsymbol{\lambda}} p_{B,\boldsymbol{\lambda}} q_{B,\boldsymbol{\lambda}} e_{0,E}^{-1}(\boldsymbol{\lambda}),$$

which by the Cauchy–Schwarz inequality and (66) is (uniformly) bounded in probability. Adding in the diagonal term in the double sum, we see that the quantity above is (uniformly) within $o_\tau(1)$ of

$$g_E^{UV}(\boldsymbol{\lambda}) - pf^{-1}g_E^U(\boldsymbol{\lambda})g_E^V(\boldsymbol{\lambda})e_{0,E}^{-1}(\boldsymbol{\lambda}).$$

Part (e) now follows from (c) and (d).

For part (f), note the given expression has conditional mean zero given $(E,\eta)$, and apply part (e) with $v_j = u_j$.

For (g), let for $\varepsilon \in (0,1)$, $\Gamma_N = \Gamma_N^{\varepsilon\uparrow} + \Gamma_N^{\varepsilon\downarrow}$, where

$$\Gamma_N^{\varepsilon\uparrow} = \frac{1}{N}\sum_{w_j \geq \varepsilon}(u_j - \overline{u}_N)^{\otimes 2} \quad \text{and} \quad \Gamma_N^{\varepsilon\downarrow} = \frac{1}{N}\sum_{w_j < \varepsilon}(u_j - \overline{u}_N)^{\otimes 2}$$

and similarly define

$$u_N^{\varepsilon\uparrow} = \frac{1}{N}\sum_{w_j \geq \varepsilon}u_j \quad \text{and} \quad u_N^{\varepsilon\downarrow} = \frac{1}{N}\sum_{w_j < \varepsilon}u_j.$$

Applying Hölder's inequality,

$$|\Gamma_N^{\varepsilon\downarrow}| \leq \left(\frac{1}{N}\sum_{j\in\mathcal{R}}|u_j - \overline{u}_N|^3\right)^{2/3}\left(\frac{1}{N}\sum_{j\in\mathcal{R}}\mathbf{1}_{(w_j<\varepsilon)}\right)^{1/3},$$

since $|u_j|^3$ is asymptotically dominated, $\Gamma_N^{\varepsilon\downarrow}$ can be made arbitrarily small by choice of $\varepsilon$. Hence, letting $\gamma$ be the value of the lim inf in (70), there exists $\varepsilon > 0$ such that

(76) $\quad \liminf\limits_{N\to\infty} \inf\limits_{|\mathbf{a}|=1} \mathbf{a}^\mathsf{T}\Gamma_N^{\varepsilon\uparrow}\mathbf{a} > 2\gamma/3 \quad \text{and} \quad |\overline{u}_N^{\varepsilon\downarrow}|^2 < \gamma/3 \qquad \text{for all } N.$

With $\geq$ the standard partial ordering on positive definite matrices, for any $\mathcal{S}$,

$$\sum_{j\in\mathcal{S}}(u_j - v)^{\otimes 2} \geq \sum_{j\in\mathcal{S}}(u_j - \overline{u}_\mathcal{S})^{\otimes 2} \qquad \text{for } \overline{u}_\mathcal{S} = \frac{1}{|\mathcal{S}|}\sum_{j\in\mathcal{S}}u_j \text{ and all } v \in \mathbf{R}^d,$$

so for this $\varepsilon$,

$$\Gamma_{N,w} \geq \sum_{w_j\geq\varepsilon}(u_j - \overline{u}_{N,w})^{\otimes 2}\frac{w_j}{\sum_j w_j} \geq \frac{\varepsilon}{\overline{\overline{w}}}\frac{1}{N}\sum_{w_j\geq\varepsilon}(u_j - \overline{u}_{N,w})^{\otimes 2}$$

$$\geq \frac{\varepsilon}{\overline{\overline{w}}}\frac{1}{N}\sum_{w_j\geq\varepsilon}(u_j - \overline{u}_N^{\varepsilon\uparrow})^{\otimes 2} \geq \frac{\varepsilon}{\overline{\overline{w}}}(\Gamma_N^{\varepsilon\uparrow} - (\overline{u}_N^{\varepsilon\downarrow})^{\otimes 2}) > 0$$

by (76). Since the weights $w_j = q_{j,\lambda_f}p_{j,\lambda_0}$ satisfy the given conditions, $\Gamma > 0$ by (62) of Condition 5.2. $\square$

THEOREM 5.1. *Consider a study base $\mathcal{R}$ of $N$ individuals with disease probability given by the proportional odds model* (1) *and a case control sampling design giving rise to the likelihood* (4). *If Conditions* 1.1 *and* 5.1–5.3 *are satisfied, there exists a consistent and asymptotically normal sequence $\hat{\boldsymbol{\beta}}_N$ of roots of the likelihood equation $\mathcal{L}_N(\boldsymbol{\beta}) = 0$; with $\Sigma$ as in* (71),

$$\hat{\boldsymbol{\beta}}_N \xrightarrow{p} \boldsymbol{\beta}_0 \quad \text{and} \quad \sqrt{N}(\hat{\boldsymbol{\beta}}_N - \boldsymbol{\beta}_0) \xrightarrow{d} \mathcal{N}(0, \Sigma^{-1}).$$

PROOF. We follow Theorem VI.I.I in Andersen, Borgan, Gill and Keiding (1993) from Billingsley (1961). For consistency it suffices to show that as $N \to \infty$,

(77) $$N^{-1}\mathcal{U}(\boldsymbol{\beta}_0) \xrightarrow{p} 0, \qquad N^{-1}\mathcal{I}(\boldsymbol{\beta}_0) \xrightarrow{p} \Sigma,$$

and, with $R(\boldsymbol{\beta}) = \partial \mathcal{I}(\boldsymbol{\beta})/\partial \boldsymbol{\beta}$, that there is a finite constant $K$ such that for some neighborhood $\mathcal{B}_0 \subset \mathcal{B}$ of $\boldsymbol{\beta}_0$,

(78) $$\lim_{N \to \infty} \mathbb{P}(|N^{-1} R(\boldsymbol{\beta})| \leq K \text{ for all } \boldsymbol{\beta} \in \mathcal{B}_0) = 1.$$

The first claim in (77) and that $N^{-1}$ times (7) tends to zero in probability follow from Lemma 5.1, part (f), Condition 5.2 and the fact that $N/|E| \xrightarrow{p} f/p$. The second claim in (77) now follows from (6) and Lemma 5.1, part (e).

Turning to (78), write, for example, $\mathbf{Z_d}$ for $\sum_{A \in \mathbf{d}} \mathbf{Z}_A$, so by (5),

(79) $$R(\boldsymbol{\beta}) = (-\mathbf{Z}''_\mathbf{D} + \mathbb{E}_{E,\eta}(\mathbf{Z}''_\mathbf{D})) + 2\operatorname{Cov}_{E,\eta}(\mathbf{Z_D}, \mathbf{Z}'_\mathbf{D})$$
$$+ \operatorname{Cov}_{E,\eta}(\mathbf{Z}'_\mathbf{D}, \mathbf{Z_D}) + \mathbb{E}_{E,\eta}\mathcal{U}(\boldsymbol{\beta})^{\otimes 3}.$$

Divided by $N$, the term inside the first parentheses tends to zero uniformly in probability over $\mathcal{B}_0$ by Lemma 5.1, part (f) and Condition 5.2. The covariances are uniformly bounded in probability upon division by $N$ by Lemma 5.1, part (e).

Last, the final term (79) over $|E|$ expands to terms of three types. For the diagonal,

$$\left| \frac{1}{|E|} \sum_{A \in E} \mathbf{Z}_A^{\otimes 3} \mathbb{E}_{E,\eta}(I_A - p_A)^3 \right| \leq \frac{1}{|E|} \sum_{A \in E} |\mathbf{Z}_A|^3,$$

for the double sums of the following form apply Corollary 4.1 to see that

$$\left| \frac{1}{|E|} \sum_{|\{A,B\}|=2} \mathbf{Z}_A^{\otimes 2} \otimes \mathbf{Z}_B \mathbb{E}_{E,\eta}(I_A - p_A)^2(I_B - p_B) \right|$$
$$\leq \frac{O_\tau(1)}{|E|^2} \sum_{|\{A,B\}|=2} |\mathbf{Z}_A|^2 |\mathbf{Z}_B|$$
$$\leq O_\tau(1) \left( \frac{1}{|E|} \sum_{A \in E} |\mathbf{Z}_A|^2 \right) \left( \frac{1}{|E|} \sum_{B \in E} |\mathbf{Z}_B| \right),$$

and for the triple sums, by Corollary 4.1,

$$\left| \frac{1}{|E|} \sum_{|\{A,B,C\}|=3} \mathbf{Z}_A \otimes \mathbf{Z}_B \otimes \mathbf{Z}_C \, \mathbb{E}_{E,\eta}(I_A - p_A)(I_B - p_B)(I_C - p_C) \right|$$

$$\leq \frac{O_\tau(1)}{|E|^3} \sum_{|\{A,B,C\}|=3} |\mathbf{Z}_A| |\mathbf{Z}_B| |\mathbf{Z}_C| \leq O_\tau(1) \left( \frac{1}{|E|} \sum_{A \in E} |\mathbf{Z}_A| \right)^3.$$

These terms are uniformly bounded over $\mathcal{B}_0$ by Condition 5.2, giving the existence of the required $K$ in (78) and completing the proof of consistency.

By the Cramér–Wold device, to prove the asymptotic normality claim, it suffices to show

(80) $\qquad \dfrac{1}{\sqrt{N}} \mathbf{b}' \mathcal{U}(\boldsymbol{\beta}_0) \xrightarrow{d} \mathcal{N}(0, \mathbf{b}' \Sigma \mathbf{b}) \qquad$ for all nonzero $\mathbf{b} \in \mathbb{R}^p$.

For $\varepsilon > 0$, define

$$\sigma_E^2 = \frac{f|E|}{Np} \mathbf{b}' \left( e_{2,E}(\boldsymbol{\lambda}) - \frac{e_{1,E}^{\otimes 2}(\boldsymbol{\lambda})}{e_{0,E}(\boldsymbol{\lambda})} \right) \mathbf{b},$$

$$G_{\varepsilon,E} = \left\{ A \in E : \left| \mathbf{b}' \mathbf{Z}_A - \frac{\mathbf{b}' e_{1,E}(\boldsymbol{\lambda})}{e_{0,E}(\boldsymbol{\lambda})} \right| > \varepsilon \sigma_E \sqrt{N} \right\},$$

$$L_{\varepsilon,E} = \frac{1}{N \sigma_E^2} \sum_{A \in G_{\varepsilon,E}} \left( \mathbf{b}' \mathbf{Z}_A - \frac{\mathbf{b}' e_{1,E}(\boldsymbol{\lambda})}{e_{0,E}(\boldsymbol{\lambda})} \right)^2 p_{A,\boldsymbol{\lambda}} q_{A,\boldsymbol{\lambda}}$$

and

$$\varepsilon_E^* = \inf\{\varepsilon : L_{\varepsilon,E} \leq \varepsilon\}.$$

Hájek's (1964) CLT, with the variables $y_A$ replaced by $\mathbf{b}' \mathbf{Z}_A p_{A,\boldsymbol{\lambda}}$, gives (80) if $\varepsilon_E^* \xrightarrow{P} 0$ as $N \to \infty$. By Hölder's inequality,

$$\frac{1}{|E|} \sum_{A \in E} |\mathbf{Z}_A|^2 \mathbf{1}(|\mathbf{Z}_A| > \varepsilon |E|^{1/2})$$

$$\leq \left( \frac{1}{|E|} \sum_{A \in E} |\mathbf{Z}_A|^3 \right)^{2/3} \left( \frac{1}{|E|} \sum_{A \in E} \mathbf{1}(|\mathbf{Z}_A| > \varepsilon |E|^{1/2}) \right)^{1/3},$$

which tends to zero in probability for all $\varepsilon > 0$ by Conditions 5.2 and 5.1 and Lemma 5.1, part (a). Since $\sigma_E^2$ is of order $O(1)$ in probability, $L_{\varepsilon,E} \xrightarrow{p} 0$. □

Turning now to the unconditional logisitic likelihood, for simplicity we parameterize $\lambda = \exp(\alpha)$, let $\alpha_{E,\eta} = \log(\boldsymbol{\lambda})$ and recall that $\tilde{\boldsymbol{\beta}}_N$ maximizes (8).

THEOREM 5.2. *Under the conditions of Theorem 5.1,*

$$
\sqrt{N}\begin{pmatrix} \tilde{\alpha}_N - \alpha_{E,\eta} \\ \tilde{\boldsymbol{\beta}}_N - \boldsymbol{\beta}_0 \end{pmatrix} \xrightarrow{d} \mathcal{N}\left(0, \Upsilon^{-1} - \begin{bmatrix} e_0^{-1} & 0^{\mathsf{T}} \\ 0 & 0 \end{bmatrix}\right),
\tag{81}
$$

*where*

$$
\Upsilon = \begin{bmatrix} e_0 & e_1^{\mathsf{T}} \\ e_1 & e_2 \end{bmatrix}.
$$

PROOF. We proceed as in the proof of Theorem 5.1. Since

$$
\frac{\partial \log(1+\lambda x_A)}{\partial \alpha} = p_{\lambda,A} \quad \text{and} \quad \frac{\partial p_{\lambda,A}}{\partial \alpha} = p_{\lambda,A} q_{\lambda,A},
$$

taking first and second partial derivatives of the logarithm of (8) with respect to $(\alpha, \boldsymbol{\beta})$, the unconditional logistic score and information are given, respectively, by

$$
\widetilde{\mathcal{U}}(\lambda, \boldsymbol{\beta}) = \sum_{A \in E} (I_A - p_{A,\lambda}) \begin{bmatrix} 1 \\ \mathbf{Z}_A \end{bmatrix}
\tag{82}
$$

and

$$
\tilde{\mathcal{I}}(\lambda, \boldsymbol{\beta}) = \sum_{A \in E} \begin{bmatrix} 1 & \mathbf{Z}_A^{\mathsf{T}} \\ \mathbf{Z}_A & \mathbf{Z}_A^{\otimes 2} \end{bmatrix} p_{A,\lambda} q_{A,\lambda} - \begin{bmatrix} 0 & 0^{\mathsf{T}} \\ 0 & \sum_{A \in E}(I_A - p_{A,\lambda})\mathbf{Z}'_A \end{bmatrix}.
\tag{83}
$$

By (82) and (46),

$$
\widetilde{\mathcal{U}}(\boldsymbol{\lambda}, \boldsymbol{\beta}_0) = \begin{bmatrix} 0 \\ \mathcal{U}(\boldsymbol{\beta}_0) \end{bmatrix} + \begin{bmatrix} 0 \\ \sum_{A \in E} \mathbf{Z}_A(p_A - p_{A,\boldsymbol{\lambda}}) \end{bmatrix}.
\tag{84}
$$

By Corollary 4.1, $p_A - p_{A,\boldsymbol{\lambda}} = O_p(N^{-1})$, so by (a) of Lemma 5.1,

$$
\sum_{A \in E} \mathbf{Z}_A(p_A - p_{A,\boldsymbol{\lambda}}) = O_p(1).
\tag{85}
$$

In view of (77), $N^{-1}\widetilde{\mathcal{U}}(\boldsymbol{\lambda}, \boldsymbol{\beta}_0) \xrightarrow{p} 0$. Handling the second term in $\tilde{\mathcal{I}}(\boldsymbol{\lambda}, \boldsymbol{\beta}_0)$ in this same manner and applying (c) of Lemma 5.1 to the first,

$$
N^{-1}\tilde{\mathcal{I}}(\boldsymbol{\lambda}, \boldsymbol{\beta}_0) \xrightarrow{p} \Upsilon.
$$

By bounding $p_{j,\boldsymbol{\lambda}} q_{j,\boldsymbol{\lambda}}$ below and following a similar but simpler argument as in (g) of Lemma 5.1, we have that $\Upsilon > 0$ by (62) of Condition 5.2.

Next we consider the remainder term. Writing $\mathbf{y}^{\mathsf{T}} = (1, \mathbf{z}^{\mathsf{T}})$ and $\boldsymbol{\gamma}^{\mathsf{T}} = (\alpha, \boldsymbol{\beta}^{\mathsf{T}})$, taking the derivative of $\tilde{\mathcal{I}}$ with respect to $\boldsymbol{\gamma}$ yields

$$
\widetilde{\mathcal{R}}(\boldsymbol{\gamma}) = \sum_{A \in E} Y_A^{\otimes 3}(p_{A,\lambda} q_{A,\lambda}^2 - p_{A,\lambda}^2 q_{A,\lambda}) + \sum_{A \in E}(Y'_A \otimes Y_A + Y_A \otimes Y'_A) p_{A,\lambda} q_{A,\lambda}
$$
$$
+ \sum_{A \in E}(I_A - p_{A,\lambda})Y''_A - p_{A,\lambda} q_{A,\lambda} Y_A \otimes Y'_A,
$$

of which all terms, once divided by $N^{-1}$, are uniformly asymptotically dominated by Condition 5.2.

By (84), (85), (80) and Slutsky's lemma,

$$(86) \qquad N^{-1/2}\widetilde{\mathcal{U}}(\boldsymbol{\lambda}, \boldsymbol{\beta}_0) \xrightarrow{d} \mathcal{N}(0, V) \qquad \text{where } V = \begin{bmatrix} 0 & 0^{\mathsf{T}} \\ 0 & \Sigma \end{bmatrix}.$$

The proof is completed by applying the well-known partitioned matrix inverse formula,

$$\Upsilon^{-1} = \begin{bmatrix} e_0^{-1} + e_0^{-1} e_1^{\mathsf{T}} \Sigma^{-1} e_1 e_0^{-1} & -e_0^{-1} e_1^{\mathsf{T}} \Sigma^{-1} \\ -\Sigma^{-1} e_1 e_0^{-1} & \Sigma^{-1} \end{bmatrix},$$

and observing

$$\Upsilon^{-1} V \Upsilon^{-1} = \begin{bmatrix} e_0^{-1} e_1^{\mathsf{T}} \Sigma^{-1} e_1 e_0^{-1} & -e_0^{-1} e_1^{\mathsf{T}} \Sigma^{-1} \\ -\Sigma^{-1} e_1 e_0^{-1} & \Sigma^{-1} \end{bmatrix} = \Upsilon^{-1} - \begin{bmatrix} e_0^{-1} & 0^{\mathsf{T}} \\ 0 & 0 \end{bmatrix}. \qquad \square$$

We note from Theorems 5.1 and 5.2 the conditional and unconditional logistic maximum likelihood estimators of the odds ratio parameter $\boldsymbol{\beta}$ have the same asymptotic distribution since $(\Upsilon^{-1})_{\beta,\beta} = \Sigma^{-1}$.

The following specialization of Theorems 5.1 and 5.2 is a direct consequence of the law of large numbers.

THEOREM 5.3. *Let $\mathbf{Z}_j, j \in \mathcal{R}$, be i.i.d. replicates of $\mathbf{Z}$. Then the conclusions of Theorems 5.1 and 5.2 hold when Conditions 1.1, 5.1 and 5.3 are satisfied, $E|\mathbf{Z}_j|^4 < \infty$, there exists an integrable random variable which bounds $|\mathbf{Z}_j|^3, |\mathbf{Z}'_j|^2$ and $|\mathbf{Z}''_j|^2$ in a neighborhood $\mathcal{B}_0 \subset \mathcal{B}$ of $\boldsymbol{\beta}_0$, and $\text{Var}(\mathbf{Z})$ is positive definite.*

When $\mathbf{Z}_j$ in the study base are independent with common distribution $\mathbf{Z}_j \stackrel{d}{=} \mathbf{Z}$, where $\mathbf{Z}$ has distribution function $G$, the case-control set $(E, \eta)$ consists of $\eta$ and $|E| - \eta$ covariates with distribution functions $G_1, G_0$, respectively, where

$$dG_i(\mathbf{z}) = \frac{p^i(\mathbf{z})(1-p(\mathbf{z}))^{1-i}}{\mathbb{E} p^i(\mathbf{Z})(1-p(\mathbf{Z}))^{1-i}} dG(\mathbf{z}), \qquad i = 0, 1,$$

with $p(\mathbf{z})$ as in (1). Then $p = \mathbb{E} p(\mathbf{Z})$, and the asymptotic distribution of $\mathbf{Z}_A$ in the case-control study is therefore given by $G_f$, where

$$dG_f(\mathbf{z}) = f \, dG_1(\mathbf{z}) + (1-f) \, dG_0(\mathbf{z}) = \frac{1-f}{1-p} \left( \frac{1 + \lambda_f \, x(\mathbf{z})}{1 + \lambda_0 \, x(\mathbf{z})} \right) dG(\mathbf{z}),$$

and the functions $h_E(\lambda)$ and $e_{k,E}(\lambda)$ converge uniformly in probability, respectively, to $h(\lambda) = \mathbb{E}_f[p_{j,\lambda}]$ and $e_k(\lambda) = \mathbb{E}_f[\mathbf{Z}_j^{\otimes k} p_{j,\lambda} q_{j,\lambda}]$.

## 6. Discussion.

**Local central limit theorem expansion for the Poisson–Binomial distribution.** The distribution of the sum of independent Bernoulli random variables with differing probabilities of success has no simple form. Theorem 2.1 gives an expansion, with rates, to any desired accuracy.

**Rejective sampling: inclusion and correlations.** The probability that an individual is included in a simple random sample has a simple form. Theorem 3.1, which gives an expansion for the probability of inclusion in a rejective sample, shows how special the equally weighted simple random sampling special case is.

Additionally, the decay rate of the high-order correlations for inclusion in a rejective sample (9) has not been previously studied, even for simple random sampling. Theorem 4.1 shows that (with $|H| = k$) the $k$th order correlations decay at the rate $|E|^{-(k+k \bmod 2)/2}$, that is, the odd correlations decay at the same rate as the next even one. In the case of simple random sampling, we have conjectured in Section 1.2 the values of the limiting constants.

**Sampling designs.** Table 1 is a list of control sampling methods most commonly used in unmatched case-control studies. The designs are classified as "case-control" type when sampling is done directly from the controls in the study base, and as "case-base" type when the sampling is from the study base without regard to case-control status. Each can be sub-classified according to whether the sampling is by simple random sampling without replacement or by independent Bernoulli trials, and whether the number of subjects to be sampled is determined by the "observed" $|\mathbf{D}|$ or "expected" $Np$ number of cases. Each design satisfies Condition 5.3. The fourth column in the Table 1 provides the parameters for the chosen sampling design that yield asymptotic case-proportion $f$. Thus, under the stated conditions on the covariates in the study base, Theorems 5.1 and 5.2 apply for each design.

TABLE 1
*Examples of sampling methods that satisfy Condition 5.3 with the parameters to yield case-proportion $f$ in the case-control set*

| Design[a] type | Sampling[b] method | Observed[c]/ expected | Sampling method to yield $f$ |
|---|---|---|---|
| C/C | SRS | Obs | Exactly $|\mathbf{D}|(1-f)/f$ controls |
| C/C | SRS | Exp | Exactly $Np(1-f)/f$ controls |
| C/C | BT | Obs | Sample controls with prob $\frac{1-f}{f}\frac{|\mathbf{D}|}{n-|\mathbf{D}|}$ |
| C/C | BT | Exp | Sample controls with prob $\frac{1-f}{f}\frac{p}{1-p}$ |
| CB | SRS | Obs | Exactly $N\frac{1-f}{f}|\mathbf{D}|/(N-|\mathbf{D}|)$ from study base |
| CB | SRS | Exp | Exactly $N\frac{1-f}{f}p/(1-p)$ from study base |
| CB | BT | Obs | Sample with prob $\frac{1-f}{f}|\mathbf{D}|/(N-|\mathbf{D}|)$ from study base |
| CB | BT | Exp | Sample with prob $\frac{1-f}{f}p/(1-p)$ from study base |

[a]C/C—case-control, CB—case-base

[b]SRS—simple random sampling, BT—Bernoulli trials

[c]Observed or expected number of cases

**Conditional and unconditional logistic regression.** Theorems 5.1 and 5.2 provide the asymptotics of the conditional and unconditional logistic likelihood estimators of the odds ratio parameter under very broad conditions. The asymptotics for the conditional estimator for this wide variety of sampling schemes are new; see Table 1. Those for the unconditional estimator extend its validity to a much wider range of applications.

Under Conditions 1.1 and 5.1–5.3, these two estimators have the same asymptotic distribution. Thus, from a statistical efficiency standpoint, either may be used. Generally, permutation likelihoods are computationally quite intensive, with complexity increasing exponentially with sample size [Liang and Qin (2000)]. However, exploiting the simplifications possible with a dichotomous outcome, a recursive algorithm for the conditional logistic likelihood reduces the order of computation to linear in $\eta$ [Cox (1972) and Gail, Lubin and Rubinstein (1981)], the same order as for the unconditional logistic likelihood. This algorithm has been implemented in a number of computer software packages. Since the unconditional estimator is biased when the number of cases is small [Breslow and Day (1980)], the conditional estimator may be preferred in situations where the case-control study consists of multiple case-control sets, some with small numbers of cases.

**Comparison to the analysis of individually matched case-control studies.** In earlier work, we studied the asymptotic behavior of conditional logistic (partial likelihood) estimators of the rate ratio from individually matched (nested) case-control data [Goldstein and Langholz (1992) and Borgan, Goldstein and Langholz (1995)].

In the individually matched case-control setting, the within case-control set variability is constant with sample size and the asymptotics are driven by the increasing number of case-control sets. The situation for the unmatched case-control setting that we studied here is very different. There is a single (or a fixed number, see Extensions below) case-control set, and the number of cases in the set increases with sample size. Consequently, a very different set of analytic techniques is required for individually matched and unmatched case-control study designs.

**Comparison to the retrospective model.** It is of interest to compare our development of the asymptotic theory of the unconditional logistic estimator to that developed under the retrospective model by Prentice and Pyke (1979). In contrast to our results, which only require asymptotic stability of the covariates, the asymptotic theory developed under the retrospective model assumes that the $\mathbf{Z}_A$ are random variables with realizations that are i.i.d. conditional on the failure indicator $I_A$. As the antihypertensive drug-MI study example in Section 1 illustrates, the identical distribution assumption may not hold in practice.

Furthermore, we note that the retrospective model is actually semiparametric, the unknown parameters being $(G_0, \boldsymbol{\beta})$, the control covariate distribution and the odds ratio parameter. Hence, efficiency questions regarding this model must be addressed by considering $G_0$ as an infinite-dimensional nuisance parameter. On the other hand, the nested case-control model considered here is parametric, leaving such questions amenable to simpler analysis.

Interestingly, the derivation of the asymptotic theory in Prentice and Pyke (1979) is quite different from the one given here. In particular, up to the scaling factor of $f/p$ which appears here, the asymptotic information $\Upsilon$ is the same for both models but the asymptotic variance of the score under the retrospective model is $\Upsilon - (f^{-1} + (1-f)^{-1})[e_0 \ e_1^\mathsf{T}]^\mathsf{T}[e_0 \ e_1^\mathsf{T}]$, compared to $V$ in (86). In spite of this difference, the asymptotic distribution of the estimator $\tilde{\boldsymbol{\beta}}$ obtained using the unconditional logistic likelihood is the same under both models. The same is almost true for $\tilde{\alpha}$, except that $e_0^{-1}$ in the nested case-control model variance (81) is replaced by $f^{-1} + (1-f)^{-1}$ in the retrospective model variance [Prentice and Pyke (1979), page 408], the difference being explained by the choice of centering values, which here is $\alpha_{E,\eta}$, and in Prentice and Pyke (1979) is $\delta$. Noting that $(f^{-1} + (1-f)^{-1})^{-1} = f(1-f) = \mathbb{E}_f(p_{A,\lambda_f})\mathbb{E}_f(q_{A,\lambda_f})$ and that $e_0 f/p = \mathbb{E}_f(p_{A,\lambda_f} q_{A,\lambda_f})$, it can be shown that the nested case-control variance associated with $\tilde{\alpha}$ is smaller than its retrospective model counterpart due to the extra conditioning here on the $\mathbf{Z}_A$.

**Efficiency.** The maximum unconditional logistic likelihood estimator has been shown to be efficient under the retrospective model [Breslow, Robins and Wellner (2000)]. These authors assume that $(I_A, \mathbf{Z}_A)$ are i.i.d., a somewhat more restricted setting than that considered by Prentice and Pyke (1979). An open question is under what conditions are $\hat{\boldsymbol{\beta}}$ and $\tilde{\boldsymbol{\beta}}$ efficient for all designs that satisfy Condition 5.3 under the nested case-control model. It would seem that the number of cases $\eta$ has no information about $\boldsymbol{\beta}_0$ so that the likelihood $\mathbb{P}_{\lambda,\beta}(\mathbf{D}|E,\eta)$ conditioning additionally on $\eta$ should not result in loss of information relative to the likelihood $\mathbb{P}_{\lambda,\beta}(\mathbf{D}|E)$. The asymptotic theory for estimators based on $\mathbb{P}_{\lambda,\beta}(\mathbf{D}|E)$ has not yet been developed (indeed, the results in this paper are a relevant step to develop such theory) so that it is not possible to compare. However, we show that the asymptotic variances of the odds ratio estimators based on $\mathbb{P}_{\lambda,\beta}(\mathbf{D}|E,\eta)$ and $\mathbb{P}_{\lambda,\beta}(\mathbf{D}|E)$ are equal in the following three important special cases.

*Simple random sampling of controls.* This class of designs, where a fixed number of controls is sampled from the study base, includes frequency matching and sampling a fixed number of controls proportional to the "expected" number of cases (i.e., the C/C-SRS entries in Table 1). For these designs, the number of cases is a function of the number in the case-control set so that $\mathbb{P}_{\lambda,\beta}(\mathbf{D}|E) = \mathbb{P}_{\lambda,\beta}(\mathbf{D}|E,\eta)$.

*Full study base.* Condition 5.3 clearly holds with $f = p$, so noting that the full, efficient likelihood $\mathbb{P}_{\lambda,\beta}(\mathbf{D}|\mathcal{R})$ has the form of an unconditional logistic likelihood, by Theorem 5.2 and 5.1 both the conditional and unconditional likelihood are efficient for $\boldsymbol{\beta}$, and in particular have the same asymptotic variance.

*Independent Bernoulli trials sampling of controls with probability $\rho$.* Under this independent control sampling design, Condition 5.3 holds with $f = p/(p + (1-p)\rho)$ and $\mathbb{P}_{\lambda,\beta}(\mathbf{D}|E)$ has the form of an unconditional logistic likelihood, and the desired conclusion follows as for the full study base.

**Extensions.** The extension to sampling controls from each of a fixed number of large strata is straightforward. Consider a failure probability model given by (1), with baseline odds parameters $\lambda_s$ for individuals in stratum $s$, and control selection independent between strata. For each $s$, let $\boldsymbol{\lambda}_s$ be the solution to

$$h_{E_s}(\lambda) = \frac{1}{|E_s|} \sum_{A \in E_s} p_{A,\lambda} = \frac{\eta_s}{|E_s|} \quad \text{and} \quad e_{k,E_s}(\lambda) = \frac{1}{|E_s|} \sum_{A \in E_s} \mathbf{Z}_A^{\otimes k} p_{A,\lambda} q_{A,\lambda},$$

where $E_s$ and $\eta_s$ are the case-control set and the number of cases from stratum $s$, respectively. Suppose the limiting fractions $\gamma_s$ of subjects in stratum $s$ exist and are positive, and that Conditions 5.1–5.3 are satisfied by all strata. Then the conclusions of Theorems 5.1 and 5.2 hold with $\Sigma = \sum_s \gamma_s \Sigma_s$, where $\Sigma_s$ is the stratum $s$ contribution to the score of form (71).

Usually disease is rare and efforts are made to enroll all cases into a case-control study. The reasons that cases are not enrolled may depend on a variety of factors, including the death of the patient or physician refusal. If nonenrollment can be modeled as i.i.d. Bernoulli($\rho_{\text{case}}$) events, then the

REJECTIVE SAMPLING ASYMPTOTICS 43

theory can easily be extended to accommodate such case selection. Specifically, the probability of sampling the case is absorbed into the baseline by replacing $\lambda_0$ in (1) by $\lambda_0 \rho_{\text{case}}$; the theory proceeds without further change.

We have used the "observed information," $-\partial \mathcal{U}(\boldsymbol{\beta})/\partial \boldsymbol{\beta}$, in our analysis. In data analysis, it is more common to use the "expected information," which is the conditional expectation over case-occurrence of the information $\mathcal{I}_{E,\eta}$ and $\tilde{\mathcal{I}}_E$ for the conditional and unconditional likelihoods, given in (6), (7) and (83), respectively [Thomas (1981)]. Because taking this expectation eliminates the term (7) in $\mathcal{I}_{E,\eta}$ and a corresponding term in $\tilde{\mathcal{I}}_E$ that was asymptotically negligible, it is immediate that the "expected information" is a consistent estimator of the asymptotic information.

**Further work.** That $\lambda_f = \lambda_0/\rho_f$ suggests that $\lambda_0$ can be estimated using the unconditional logistic likelihood when the number of subjects in the study base (and thus the proportion of cases) is known. This has been done under (essentially) Bernoulli trials by Weinberg and Wacholder (1993), and under independent simple random sampling of (cases and) controls by Scott and Wild (1986) and Breslow and Cain (1988), but further work is needed to accommodate general Condition 5.3 sampling. In particular, there is nonnegligible variability in the difference $\boldsymbol{\lambda} - \lambda_f$ that depends on the sampling design, and which needs to be accounted for in the estimation of $\lambda_0$.

It is of interest to know when the techniques used here can be generalized to accommodate other forms of conditioning on information $\mathcal{S}$, as in likelihood (3). The particular case of no conditioning, $\mathcal{S} = \varnothing$, represents a "full likelihood" under the nested case-control model. The difficulty is finding an analog to the independent product measure $T_\lambda$ in Lemma 3.5.

R. Arratia
L. Goldstein
Department of Mathematics
University of Southern California
Kaprielian Hall, Room 108
3620 Vermont Avenue
Los Angeles, California 90089-2532
USA
e-mail: larry@math.usc.edu

B. Langholz
Department of Preventive Medicine
Keck School of Medicine
1540 Alcazar St., CHP-220
Los Angeles, California 90089-9011
USA
e-mail: langholz@usc.edu